\documentclass[amssymb,11pt]{amsart}
\usepackage{amscd,amssymb,euscript,amsthm}
\theoremstyle{plain}
\newtheorem{thm}{Theorem}[section] 
\newtheorem{cor}[thm]{Corollary}
\newtheorem{prop}[thm]{Proposition}
\newtheorem{lem}[thm]{Lemma}

\theoremstyle{definition}
\newtheorem{defn}[thm]{Definition}
\theoremstyle{remark}
\newtheorem{rem}[thm]{Remark}

\numberwithin{equation}{section}
\newcommand{\tr}{\operatorname{trace}}

\newcommand{\co}{\operatorname{conv}}

\def\<{\left<}
\def\>{\right>}
\def\cstar{$C^*$-algebra}
\begin{document}
\title[Maximal vectors]{Maximal vectors in Hilbert space\\
and quantum entanglement}
\author{William Arveson}
%
%
%
%
\date{12 May, 2008}

\begin{abstract}  Let $V$ be a norm-closed subset of the unit sphere 
of a Hilbert space $H$ that is stable under multiplication 
by scalars of absolute value $1$.    
A {\em maximal vector} (for $V$) is a unit vector $\xi\in H$ whose distance 
to $V$ is maximum 
$$
d(\xi,V)=\sup_{\|\eta\|=1}d(\eta,V),  
$$
$d(\xi,V)$ denoting the distance from $\xi$ to the set $V$.  
Maximal vectors generalize the {\em maximally entangled} 
unit vectors of quantum theory.  

In general, under a mild regularity hypothesis on $V$, 
there is a {\em norm}  on $H$ whose restriction to 
the unit sphere achieves its minimum 
precisely on $V$ and its 
maximum precisely on the set of maximal vectors.  This 
``entanglement-measuring norm"  is
unique.  There is a corresponding ``entanglement-measuring norm" 
on the predual of $\mathcal B(H)$ that faithfully detects entanglement 
of normal states.  

We apply these abstract results to the analysis of 
entanglement in multipartite 
tensor products $H=H_1\otimes \cdots\otimes H_N$, 
and we calculate both entanglement-measuring norms.  
In  cases for which $\dim H_N$ is 
relatively large with respect to the others,  
we describe the set of maximal 
vectors in explicit terms and show that it 
does not depend on the number of factors of the 
Hilbert space $H_1\otimes\cdots\otimes H_{N-1}$. 
\end{abstract}
\maketitle

\section{Introduction}\label{S:in}

Let $H=H_1\otimes\cdots\otimes H_N$ be a finite tensor product of separable Hilbert spaces.  
In the literature of physics and quantum information theory, 
a normal state $\rho$ of $\mathcal B(H_1\otimes\cdots\otimes H_N)$ called {\em separable} 
or {\em classically correlated} 
if it belongs to the norm closed convex set generated by product states  
$\sigma_1\otimes \cdots\otimes \sigma_N$, where $\sigma_k$ denotes a normal state of $\mathcal B(H_k)$.  
Normal states that are not separable are said to be {\em entangled}.  The notion 
of entanglement is a distinctly noncommutative phenomenon, and 
has been a fundamental theme of quantum physics since the early 
days of the subject.  It has received increased attention recently because of 
possible 
applications emerging from quantum information theory.   

In the so-called bipartite case in which $N=2$, several numerical measures of 
entanglement have been proposed that emphasize various features 
(see  \cite{HorSurvey}, \cite{hylEtAl}, 
\cite{peresSep}, \cite{wangGu}).  Despite the variety 
of proposed measures, 
only one we have seen (the projective cross norm 
introduced in \cite{rudSep}, \cite{rudEn} )
is capable of distinguishing between entangled mixed states and 
separable mixed states of bipartite tensor products. Of course, the bipartite case has 
special features because vectors 
in $H_1\otimes H_2$ can be identified 
with Hilbert-Schmidt operators from $H_1$ to $H_2$, 
thereby allowing one to access operator-theoretic invariants -- most 
notably the 
singular value list of a Hilbert-Schmidt operator -- to 
analyze vectors in $H_1\otimes H_2$.  On the other hand, 
that tool is much less effective for higher order tensor 
products,  and perhaps that explains why the higher order 
cases $N\geq 3$ are poorly understood.  For example,  there 
does not appear to be 
general agreement as to what properties a ``maximally entangled" vector should 
have in such cases; and 
in particular, there is no precise definition of the term.

In this paper we propose such a definition and 
introduce two numerical invariants (one 
for vectors and one for states) that 
faithfully detect entanglement, in a general 
mathematical setting that includes the cases of physical interest.  
We start  with a separable Hilbert space $H$ and a 
distinguished set 
$$
V\subseteq \{\xi\in H: \|\xi\|=1\}
$$ 
of unit vectors that satisfies the following two conditions: 

\begin{enumerate}
\item[$V1$:]  $\lambda\cdot V\subseteq V$, for every $\lambda\in \mathbb C$ with $|\lambda|=1$.  
\item[$V2$:] For every $\xi\in H$,  $\langle\xi,V\rangle=\{0\}\implies \xi=0$.  
\end{enumerate}
By replacing $V$ with its closure if necessary, we can and do assume that $V$ is closed 
in the norm topology of $H$.  
A normal state $\rho$ of $\mathcal B(H)$ is said to be {\em $V$-correlated} if 
for every $\epsilon>0$, there are vectors $\xi_1,\dots,\xi_n\in V$ and 
positive numbers $t_1,\dots,t_n$ with sum $1$ such that 
$$
\sup_{\|x\|\leq 1}|\rho(x)-\sum_{k=1}^nt_k\langle x\xi_k,\xi_k\rangle|\leq \epsilon.    
$$
A normal state that is not $V$-correlated is called $V$-entangled - or simply 
{\em entangled}.  
The motivating examples are those in which $H=H_1\otimes\cdots\otimes H_N$ 
is an $N$-fold tensor product of Hilbert spaces $H_k$ and 
$$
V=\{\xi_1\otimes\cdots\otimes \xi_n: \xi_k\in H_k,\ \|\xi_1\|=\cdots=\|\xi_N\|=1\}
$$
is the set of decomposable unit vectors.  In such cases the $V$-correlated states 
are the {\em separable} states, and when $H$ is finite dimensional, the 
$V$-correlated states are the simply the convex combinations of vector 
states $x\mapsto \langle x\xi,\xi\rangle$ with $\xi$ a unit vector 
of the form $\xi=\xi_1\otimes\cdots\otimes\xi_n$, \  $\xi_k\in H_k$, 
$k=1,\dots,n$.  Of course, there are many other examples that have less 
to do with physics.  

In general, given such a set $V\subseteq H$, a {\em maximal vector} is 
defined as a 
unit vector $\xi\in H$ whose distance to $V$ is 
maximum 
$$
d(\xi,V)=\sup_{\|\eta\|=1}d(\eta,V), 
$$ 
$d(\eta,V)$ denoting the distance from $\eta$ to $V$.  While it is 
not obvious from this geometric definition, it is a fact that 
in the case of bipartite tensor products $H=H_1\otimes H_2$, maximal 
vectors turn out to be exactly the ``maximally entangled" unit vectors of 
the physics literature (see (\ref{inEq2}) below).    
Sections \ref{S:dm} through \ref{S:dn} are devoted to 
an analysis of the geometric properties of maximal vectors 
in general.  We introduce a numerical invariant $r(V)$ of $V$ (the inner radius) 
and show that when $r(V)>0$, there is a uniquely determined 
``entanglement measuring norm" $\|\cdot\|^V$ on $H$ with the property 
that $\xi\in V$ iff $\|\xi\|^V=1$ and 
$\xi$ is maximal iff $\|\xi\|^V= r(V)^{-1}$ (see Proposition \ref{dmProp1} and 
Theorem \ref{dnThmA}).  

In Section \ref{S:es} 
we introduce an extended real-valued function $E(\rho)$ of 
normal states $\rho$ that takes values in the interval 
$[1,+\infty]$.  This ``entanglement" function $E$ 
is convex, 
lower semicontinuous, and faithfully detects generalized entanglement 
in the sense that $\rho$ is entangled iff $E(\rho)>1$ (Theorems \ref{esThm1} and 
\ref{fThm1}).  We also show that under the 
same regularity hypothesis on the given set $V$ of unit vectors (namely $r(V)>0$), 
$E$ is a {\em norm} equivalent to the ambient norm of $\mathcal B(H)_*\cong\mathcal L^1(H)$, 
and it achieves its 
maximum on vector states of the form $\omega(A)=\langle A\xi,\xi\rangle$, 
$A\in\mathcal B(H)$ precisely when $\xi$ is a maximal vector 
(Theorem \ref{msThm1}).  

In the third part of the paper (Sections \ref{S:vv}--\ref{S:cn}), we apply these abstract 
results to cases in which $V$ is the 
set of decomposable unit vectors in an $N$-fold tensor product $H=H_1\otimes \cdots\otimes H_N$.   
We assume that all but one of the $H_k$ are finite dimensional, arranged so that the 
dimensions $n_k=\dim H_k$ weakly increase with $k$ and $n_{N-1}<\infty$. 
We identify the vector norm $\|\cdot\|^V$ that measures entanglement as 
the greatest cross norm on the projective tensor product of Hilbert spaces 
$$
H_1\hat\otimes H_2\hat\otimes \cdots\hat\otimes H_N 
$$
in general - see Theorem \ref{vvThm1}.  Similarly, we identify the  entanglement function 
of mixed states as the 
restriction to density operators of the greatest cross norm 
of the projective tensor product of Banach spaces 
$$
\mathcal L^1(H_1)\hat\otimes\mathcal L^1(H_2)\hat\otimes\cdots\hat\otimes \mathcal L^1(H_N), 
$$ 
$\mathcal L^1(H)$ denoting the Banach space of trace class operators 
on a Hilbert space $H$ (Theorem \ref{msThm1}).  Note that 
in the bipartite case $N=2$, the latter reduces to the norm introduced 
in a more {\em ad hoc} way by Rudolph in \cite{rudSep}, \cite{rudEn}.  

We are unable to identify the maximal vectors in this generality, and 
our sharpest results  for multipartite tensor products require an additional 
hypothesis, namely that one of the spaces $H_k$ should be 
significantly larger than the others in the sense that $n_N\geq n_1\cdots n_{N-1}$.  
In every case, of course, the entanglement measuring norm $\|\cdot\|^V$ depends 
strongly on relative dimensions $n_1,\dots, n_{N}$ of the factors of 
the decomposition 
$H_1\otimes\cdots\otimes H_{N}$, because the ``shape" of its unit ball 
$\{\xi: \|\xi\|^V\leq 1\}$ depends strongly on these relative dimensions.  
What is interesting is that 
when $n_N\geq n_1\cdots n_{N-1}$, the set of maximal vectors does  not depend on that 
finer structure.  Indeed, we show that in such cases 
the maximal vectors are precisely the vectors in $H_1\otimes\cdots\otimes H_N$ 
that can be represented 
\begin{equation}\label{inEq1}
\xi=\frac{1}{\sqrt{n_1n_2\cdots n_{N-1}}}
\sum_{k=1}^{n_1\cdots n_{N-1}}
e_k \otimes f_k
\end{equation}
where $e_1,\dots, e_{n_1\cdots n_{N-1}}$ is an orthonormal 
basis for $H_1\otimes \cdots\otimes H_{N-1}$ and $f_1,\dots, f_{n_1\cdots n_{N-1}}$ 
is an orthonormal set in $H_N$ 
(see Theorem \ref{hThm1}).  The 
simplest case is $N=2$, where our hypotheses reduce to 
$n_1\leq n_2\leq \infty$ with $n_1$ finite, 
and the expression (\ref{inEq1}) becomes a familiar 
representation of ``maximally entangled" vectors of bipartite tensor products 
that is commonly found in the physics literature.  

To make the point in somewhat more physical terms, let $H$ and $K$ be finite dimensional 
Hilbert spaces with $n=\dim H\leq m=\dim K$.  The maximal vectors of 
the bipartite tensor product $H\otimes K$ are those of the form 
\begin{equation}\label{inEq2}
\xi=\frac{1}{\sqrt n}(e_1\otimes f_1+\cdots+e_n\otimes f_n)
\end{equation}
where $(e_k)$ is an orthonormal basis for $H$ and $(f_k)$ is an orthonormal set 
in $K$.  
On the other hand, if the Hilbert space $H$ represents a composite 
of several subsystems in the sense that it can be decomposed 
into a tensor product 
$H=H_1\otimes \cdots\otimes H_r$ of Hilbert spaces, then the set of  maximal 
vectors relative to the more refined decomposition $H_1\otimes\cdots\otimes H_r\otimes K$
is precisely the same set of unit vectors (\ref{inEq2}).

This unexpected stability of the set of maximal vectors is established by 
showing that the states associated 
with maximal vectors $\xi$ are  
characterized by the following requirement on their 
``marginal distributions".  
The algebra $\mathcal A=\mathcal B(H_1\otimes\cdots\otimes H_{N-1})$ 
can be viewed as a matrix algebra with tracial state $\tau$, and 
we show that a unit vector 
$\xi$ in $H_1\otimes \cdots\otimes H_N$ is maximal if and only if 
$$
\langle (A\otimes\mathbf 1_{H_N})\xi,\xi\rangle =\tau(A), \qquad A\in \mathcal A, 
$$
see Theorem \ref{scThm1}.  We do not know if there is a useful characterization 
of the marginal states of maximal vectors in the remaining cases for 
which $n_N<n_1\cdots n_{N-1}$, and that is an issue calling for further research.

Of course, it was also necessary to calculate the geometric invariant $r(V)$
for these examples,  
see Theorems \ref{ciThm1}
and \ref{ciThm2}.  A more precise and more complete summary of our main results for 
multipartite tensor products is presented at the end of the paper in Theorem \ref{srThm1} 
(also see Remark \ref{cnRem2}).  

The idea of measuring entanglement of vectors in terms of their distance to 
the decomposable vectors appears in \cite{weiGo}, and calculations are carried out 
for several examples.  While a related measure was also introduced for states, it is 
different from the one below, and there appears to be no further overlap with 
this paper.  Also see formula (22) of \cite{guhneEtAl}.  
A related operator-theoretic notion 
of entanglement for bipartite tensor products was introduced in \cite{BNT}, where 
it is 
shown essentially that a density operator that is maximally far from the separable 
ones relative to the Hilbert-Schmidt norm provides a maximal violation of the 
Bell inequalities.  Perhaps it is also relevant to point out 
that the recent paper \cite{pgEtAl} establishes unbounded violations 
of the Bell inequalities for tripartite tensor products using quite 
different methods.

This is the third of a series of papers that relate to 
entangled states 
on matrix algebras \cite{arvEnt1}, \cite{arvEnt2} .  However, while 
the results below certainly apply to 
matrix algebras, many of them also apply to the context 
of infinite dimensional Hilbert spaces.  Finally, I  thank Mary Beth Ruskai 
for calling my attention to some key results in the physics literature, and Yoram 
Gordon for helpful comments.

\part{Vectors in Hilbert spaces}

\section{Detecting membership in convex sets}\label{S:dm}

Let $H$ be a Hilbert space and let $V\subseteq \{\xi\in H: \|\xi\|=1\}$ be 
a norm-closed subset of the unit sphere of $H$ that satisfies 
V1 and V2.  Recall that since the weak closure and the norm closure of 
a convex subset of $H$ are the same, it is unambiguous to speak of 
the closed convex hull of $V$.  

In this section we show that there is 
a unique function $u: H\to [0,+\infty]$ with certain 
critical properties that determines membership in 
the closed convex hull of $V$, and more significantly for our purposes, 
{\em such a function 
determines membership in $V$ itself}.    While the proof 
of Proposition \ref{dmProp1} below involves 
some familiar ideas from convexity theory, 
it is not part of the lore of topological vector 
spaces, hence we include details.

We begin with a preliminary function $\|\cdot\|_V$ defined on $H$ by 
\begin{equation}\label{dmEq1}
\|\xi\|_V=\sup_{\eta\in V}|\langle \xi,\eta\rangle|, \qquad \xi\in H.  
\end{equation}
Axiom V2 implies that $\|\cdot\|_V$ is a norm, and 
since $V$ consists of unit vectors we have $\|\xi\|_V\leq \|\xi\|$. 
The associated unit ball $\{\xi\in H: \|\xi\|_V\leq 1\}$ 
is a closed convex subset of $H$ that contains the unit ball 
$\{\xi\in H: \|\xi\|\leq 1\}$ of $H$ because $\|\xi\|_V\leq \|\xi\|$, 
$\xi\in H$.  

Now consider the function $\|\cdot\|^V: H\to [0,+\infty]$ defined by 
\begin{equation}\label{dmEq2}
\|\xi\|^V=\sup_{\|\eta\|_V\leq 1}\Re\langle\xi,\eta\rangle=
\sup_{\|\eta\|_V\leq 1}|\langle \xi,\eta\rangle|, \qquad \xi\in H.  
\end{equation}
Since $\|\eta\|_V\leq \|\eta\|$, the right side of (\ref{dmEq2}) is at least 
$\|\xi\|$, hence 
\begin{equation}\label{dmEq3}
\|\xi\|_V\leq \|\xi\|\leq \|\xi\|^V, \qquad \xi\in H.  
\end{equation}
Significantly, it is possible for $\|\xi\|^V$ to achieve the value $+\infty$ when 
$H$ is infinite dimensional; an example is given in Proposition \ref{vvProp2} below.  

An extended real-valued function $u: H\to [0,+\infty]$ is said to be 
{\em weakly lower semicontinuous} if for every $r\in [0,+\infty)$, 
the set $\{\xi\in H: u(\xi)\leq r\}$ is closed in the weak topology of $H$.

\begin{prop}\label{dmProp1}
The  extended real-valued function 
$\|\cdot\|^V: H\to [0,+\infty]$
has the following properties:
\begin{enumerate}
\item[(i)] $\|\xi+\eta\|^V\leq \|\xi\|^V+\|\eta\|^V$, \quad $\xi,\eta\in H$.
\item[(ii)] $\|\lambda\cdot\xi\|^V=|\lambda|\cdot \|\xi\|^V$, \quad $0\neq\lambda\in\mathbb C$,\quad $\xi\in H$. 
\item[(iii)] It is weakly lower semicontinuous.  
\item[(iv)] The closed convex hull of $V$ is $\{\xi\in H: \|\xi\|^V\leq 1\}$.  
\end{enumerate}
This function is uniquely determined: 
If $u: H\to [0,+\infty]$ is any function that satisfies (ii) and (iv), 
then $u(\xi)=\|\xi\|^V$, $\xi\in H$.  
\end{prop}

The proof rests on the following result.

\begin{lem}\label{dmLem1}
Let $K$ be the closed convex hull of $V$.  Then 
\begin{equation}\label{dmEq4}
K=\{\xi\in H: \|\xi\|^V\leq 1\},   
\end{equation}
and in particular, 
\begin{equation}\label{dmEq5}
\|\xi\|_V=\sup_{\|\eta\|^V\leq 1}|\langle\xi,\eta\rangle|, \qquad \xi\in H.  
\end{equation}
\end{lem}

\begin{proof}[Proof of Lemma \ref{dmLem1}]
For the inclusion $\subseteq$ of (\ref{dmEq4}), note that if $\xi\in V$ 
and $\eta$ is any vector in $H$, then $|\langle\xi,\eta\rangle|\leq \|\eta\|_V$, 
so that 
$$
\|\xi\|^V=\sup_{\|\eta\|_V\leq 1}|\langle\xi,\eta\rangle|\leq \sup_{\|\eta\|_V\leq 1}\|\eta\|_V\leq 1.  
$$ 
For the other inclusion, a standard separation theorem implies that it is enough 
to show that for every continuous linear functional $f$ on $H$ and every $\alpha\in \mathbb R$, 
$$
\sup_{\xi\in V}\Re f(\xi)\leq \alpha\implies \sup_{\|\eta\|^V\leq 1}\Re f(\eta)\leq \alpha.  
$$
Fix such a pair $f$, $\alpha$ with $f\neq 0$.  By the Riesz lemma, there is a vector $\zeta\in H$ 
such that $f(\xi)=\langle \xi,\zeta\rangle$, $\xi\in H$, and the first inequality 
above implies  
$$
0<\|\zeta\|_V=\sup_{\xi\in V}|\langle\xi,\zeta\rangle|=\sup_{\xi\in V}\Re f(\xi)\leq \alpha.  
$$
Hence $\|\alpha^{-1}\zeta\|_V\leq 1$.  By definition of $\|\cdot\|^V$ we have 
$|\langle\eta,\alpha^{-1}\zeta\rangle|\leq \|\eta\|^V$, 
therefore $|\langle\eta,\zeta\rangle|\leq \alpha\|\eta\|^V$, and finally  
$$
\sup_{\|\eta\|^V\leq 1}\Re f(\eta)\leq \sup_{\|\eta\|^V\leq 1}|\langle \eta,\zeta\rangle|
\leq \alpha, 
$$
which is the inequality on the right of the above implication.  

To deduce the formula (\ref{dmEq5}), use (\ref{dmEq4}) to 
write 
$$
\|\xi\|_V=\sup_{\eta\in V}|\langle\xi,\eta\rangle|=\sup_{\eta\in K}|\langle \xi,\eta\rangle|
=\sup_{\|\eta\|^V\leq 1}|\langle\xi,\eta\rangle|
$$
and (\ref{dmEq5}) follows.  
\end{proof}

\begin{proof}[Proof of Proposition \ref{dmProp1}] 
Properties (i) and (ii) are obvious from the definition (\ref{dmEq2}) of $\|\cdot\|^V$, 
lower semicontinuity (iii) also follows immediately from the definition (\ref{dmEq2}),  
and property (iv) follows from Lemma \ref{dmLem1}.

 Uniqueness:  Property (iv) implies that for $\xi\in H$, 
$$
u(\xi)\leq 1\iff \|\xi\|^V\leq 1. 
$$ 
Using $u(r\cdot\xi)=r\cdot u(\xi)$ for 
$r>0$, we conclude that for every positive real number $r$ and every $\xi\in H$, one has  
$$
u(\xi)\leq r\iff \|\xi\|^V\leq r,  
$$
from which it follows that $u(\xi)=\|\xi\|^V$ whenever one of $u(\xi), \|\xi\|^V$ is finite, 
and that $u(\xi)=\|\xi\|^V=+\infty$ whenever one of $u(\xi), \|\xi\|^V$ is $+\infty$.  
Hence $u(\xi)=\|\xi\|^V$ for all $\xi\in H$.    
\end{proof}

What is more significant is that 
the function $\|\cdot\|^V$ detects membership in $V$ itself:

\begin{thm}\label{dmCor1}  The restriction of 
the function $\|\cdot\|^V$ of (\ref{dmEq2})  
to the unit sphere $\{\xi\in H: \|\xi\|=1\}$ of $H$ satisfies 
\begin{equation}\label{dmThm2}
\|\xi\|^V\geq 1, \quad \rm{and }\quad  \|\xi\|^V=1\iff \xi\in V.     
\end{equation}
\end{thm}

\begin{proof}  (\ref{dmEq3}) implies that $\|\xi\|^V\geq 1$ for all $\|\xi\|=1$.  

Let $K$ be the closed convex hull of $V$.  The 
description of $K$ given in 
(\ref{dmEq4}) and the properties (i) and (ii) of Proposition 
\ref{dmProp1} imply that the extreme points of $K$ are the vectors $\xi\in H$ satisfying 
$\|\xi\|^V=1$.  Since $V$ consists of extreme points of the unit ball of $H$, 
it consists of extreme points of $K$, hence $\|\xi\|^V=1$ for every $\xi\in V$.  

Conversely, if $\xi$ satisfies $\|\xi\|=\|\xi\|^V=1$, then the preceding 
remarks show that $\xi$ is an extreme point of $K$, so that Milman's converse 
of the Krein-Milman theorem implies that $\xi$ belongs to the weak closure 
of $V$.  But if $\xi_n$ is a sequence in $V$ that converges weakly 
to $\xi$ then 
$$
\|\xi_n-\xi\|^2=2-2\Re\langle \xi_n,\xi\rangle\to 2-2\langle \xi,\xi\rangle=0
$$
as $n\to\infty$.  We conclude that $\xi\in \overline{V}^{\,\rm{norm}}=V$.  
\end{proof}

\section{The geometric invariant $r(V)$}\label{S:cc}

In this section we introduce a numerical invariant 
of $V$ that will play 
a central role.   

\begin{defn} The {\em inner radius} $r(V)$ of $V$ is defined as the largest 
$r\geq 0$ such that $\{\xi\in H: \|\xi\|\leq r\}$ is contained 
in the closed convex hull of $V$.  
\end{defn}

Obviously, $0\leq r(V)\leq 1$.  The following result and 
its corollary imply that $r(V)>0$ when $H$ is finite dimensional.  More 
generally, they imply that whenever the inner radius is positive, 
{\em both} $\|\cdot\|_V$ and $\|\cdot\|^V$ are norms that are equivalent to the 
ambient norm of $H$.  We write $d(\xi,V)$ for the distance 
from a vector $\xi\in H$ to the set $V$, 
$d(\xi,V)=\inf\{\|\xi-\eta\|: \eta\in V\}$.

\begin{thm}\label{ccThm1}
The inner radius of $V$ is characterized by each of the following 
three formulas:  
\begin{equation}\label{ccEq2}
\inf_{\|\xi\|=1}\|\xi\|_V=r(V), 
\end{equation}
\begin{equation}\label{ccEq3}
\sup_{\|\xi\|=1}\|\xi\|^V=\frac{1}{r(V)}, 
\end{equation}
\begin{equation}\label{ccEq3.1}
\sup_{\|\xi\|=1}d(\xi,V)=\sqrt{2(1-r(V))}.  
\end{equation}
\end{thm}

\begin{proof} Let $K$ be the closed convex hull of $V$.  If $K$ contains the 
ball of radius $r$ about $0$, then for every $\xi\in H$ we have 
$$
\sup_{\eta\in V}|\langle \xi,\eta\rangle|=\sup_{\eta\in K}|\langle\xi,\eta\rangle|\geq 
\sup_{\|\eta\|\leq r}|\langle \xi,\eta\rangle|=r\cdot\|\xi\|.   
$$
Hence  
$$
\inf_{\|\xi\|=1}\|\xi\|_V=\inf_{\|\xi\|=1}\sup_{\eta\in V}|\langle\xi,\eta\rangle |\geq r,  
$$
and  $r(V)\leq \inf\{\|\xi\|_V: \|\xi\|=1\}$ follows.  
For the opposite inequality, set 
$$
r=\inf_{\|\xi\|=1}\|\xi\|_V.  
$$ 
Then for every $\xi\in H$ satisfying $\|\xi\|=1$, we have 
$$
\sup_{\|\eta\|\leq r}\Re\langle \xi,\eta\rangle=
r\cdot\sup_{\|\eta\|\leq 1}\Re\langle \xi,\eta\rangle=r\cdot\|\xi\|=r\leq 
\|\xi\|_V=\sup_{\eta\in V}\Re\langle \xi,\eta\rangle, 
$$
and after rescaling $\xi$ we obtain
$$
\sup_{\|\eta\|\leq r}\Re\langle \xi,\eta\rangle\leq \sup_{\eta\in V}\Re\langle \xi,\eta\rangle,  
\qquad \xi\in H.  
$$
At this point, a standard separation theorem implies that 
$\{\eta\in H:\|\eta\|\leq r\}$ is contained 
in the closed convex hull of $V$, hence $r\leq r(V)$.  

(\ref{ccEq3})  follows from 
(\ref{ccEq2}), since by definition of the norm $\|\xi\|^V$ 
\begin{align*}
\sup_{\|\xi\|=1}\|\xi\|^V&=\sup_{\|\xi\|=1, \|\eta\|_V=1}|\langle \xi,\eta\rangle|
=\sup_{\|\eta\|_V=1}\|\eta\|=\sup_{\eta\neq 0}\frac{\|\eta\|}{\|\eta\|_V}\\
&=\sup_{\|\eta\|=1}\frac{1}{\|\eta\|_V}
=\left(  \inf_{\|\eta\|=1}\|\eta\|_V\right)^{-1}=r(V)^{-1}.  
\end{align*}

To prove (\ref{ccEq3.1}), the distance $d(\xi,V)$ from $\xi$ to $V$ satisfies 
\begin{align*}
d(\xi,V)^2&=\inf_{\eta\in V}\|\xi-\eta\|^2=\inf_{\eta\in V}(2-2\Re\langle\xi,\eta\rangle)
=2-2\sup_{\eta\in V}\Re\langle\xi,\eta\rangle
\\
&=2-2\sup_{\eta\in V}|\langle\xi,\eta\rangle |=2-2\|\xi\|_V,  
\end{align*}
and (\ref{ccEq3.1}) follows after taking square roots.
\end{proof}

\begin{cor}\label{ccCor1}  If the inner radius $r(V)$ is positive, then 
$\|\cdot\|^V$ is a norm on $H$ satisfying 
$$
\|\xi\|\leq \|\xi\|^V\leq \frac{1}{r(V)}\|\xi\|, \qquad \xi\in H.  
$$
If $H$ is finite dimensional, then $r(V)>0$.  
\end{cor}

\begin{proof}  The first sentence follows from (\ref{dmEq3}) and (\ref{ccEq3}).  
If $H$ is finite dimensional, 
all norms on $H$ are equivalent, and $r(V)>0$ follows from 
(\ref{ccEq2}).  
\end{proof}

\begin{cor}\label{ccCor2}  In general, for any closed set 
$V$ of unit vectors that satisfies V1, the 
following five assertions about $V$ are equivalent:
\begin{enumerate}
\item[(i)] The closed convex hull of $V$ has nonempty interior.  
\item[(ii)] The inner radius of $V$ is positive.  
\item[(iii)] The seminorm $\|\cdot\|_V$ is equivalent to the ambient norm of $H$.  
\item[(iv)] The function $\|\cdot\|^V$ is a norm equivalent to the ambient norm of $H$.  
\item[(v)] The function $d(\cdot,V)$ is bounded away from $\sqrt 2$ on 
the unit sphere: 
$$
\sup_{\|\xi\|=1}d(\xi,V)<\sqrt 2.  
$$
\end{enumerate}
\end{cor}

\begin{proof}
The equivalences (ii)$\iff$(iii)$\iff$(iv)$\iff$(v) are immediate consequences 
of the formulas of Theorem \ref{ccThm1}.  Since the implication 
(ii)$\implies$(i) is trivial, it suffices to prove (i)$\implies$(ii).  

For that, let $U$ be a nonempty open set contained in the closed convex hull 
$K$ of $V$.  The vector difference $U-U$ is an 
open neighborhood of $0$ that is contained in $K-K$.  By 
axiom V1, $K-K$ is contained in $2\cdot K$, so that 
$2^{-1}\cdot(U-U)$ is a subset of $K$ that contains an open 
ball about $0$.   
\end{proof}

\section{Maximal vectors}\label{S:dn}

Throughout this section,  $V$ will denote a norm-closed subset of the unit sphere of a Hilbert space $H$ 
that satisfies V1 and V2.  For every unit vector $\xi\in H$, the distance 
from $\xi$ to $V$ satisfies $0\leq d(\xi,V)\leq \sqrt 2$; and since $V$ is norm-closed, 
one has $d(\xi,V)=0$ iff $\xi\in V$.  

\begin{defn}\label{dnDef1}
By a {\em maximal vector} we mean a vector $\xi\in H$ satisfying 
$\|\xi\|=1$ and 
$$
d(\xi,V)=\sup_{\|\eta\|=1}d(\eta,V).  
$$
\end{defn}

When $H$ is finite dimensional, an obvious compactness argument shows that maximal vectors exist;  
and they exist for significant infinite 
dimensional examples as well (see Sections \ref{S:vv}--\ref{S:sr}).  
Maximal vectors will play a central 
role throughout the remainder of this paper.  In this section we 
show that whenever $r(V)>0$,  
the restriction of the function $\|\cdot\|^V$ 
to the unit sphere of $H$ 
detects {\em maximality} as well as membership in $V$.  
Indeed, in Theorem \ref{ccThm1} we calculated 
the minimum of $\|\cdot\|_V$ 
and the maximum of $\|\cdot\|^V$ over the unit sphere of $H$.  What 
is notable is that when either of the two extremal values is achieved 
at some unit vector $\xi$ then they are both achieved at $\xi$;  
and that such vectors $\xi$ are precisely the 
maximal vectors.

\begin{thm}\label{dnThmA} 
If $r(V)>0$, then for every unit vector $\xi\in H$, the following 
three assertions are equivalent: 
\begin{enumerate}
\item[(i)] $\|\xi\|_V=r(V)$ is minimum.  
\item[(ii)] $\|\xi\|^V=r(V)^{-1}$ is maximum.  
\item[(iii)] $d(\xi,V)=\sqrt{2(1-r(V))}$ is maximum.   
\end{enumerate}
\end{thm}

\begin{proof} 
Choose a unit vector $\xi$.  
We will prove the implications (i)$\iff$(iii), (i)$\implies$(ii) and (ii)$\implies$(i).  

(i)$\iff$(iii): 
Theorem \ref{ccThm1} implies that 
the minimum value of $\|\xi\|_V$ is $r(V)$,  the maximum 
value of $d(\xi,V)$ is given by (iii), and that 
$d(\xi,V)$ is maximized at $\xi$ iff $\|\xi\|_V$ is minimized at $\xi$.

(i)$\implies$(ii):  If $\|\xi\|_V=r(V)$ then $\|r(V)^{-1}\xi\|_V=1$, so 
that 
$$
\|\xi\|^V=\sup_{\|\eta\|_V\leq 1}|\langle \xi,\eta\rangle|
\geq |\langle \xi,r(V)^{-1}\xi\rangle|=
\frac{1}{r(V)}.   
$$
Since (\ref{ccEq3}) 
implies $\|\xi\|^V\leq r(V)^{-1}$, we conclude that $\|\xi\|^V=r(V)^{-1}$.  

(ii)$\implies$(i): Assuming (ii), we have  
\begin{align*}
r(V)^{-1}&=\|\xi\|^V=\sup_{\|\eta\|_V\leq 1}|\langle\xi,\eta\rangle|
=\sup_{\|\eta\|_V=1}|\langle \xi,\eta\rangle|
\\
&=\sup_{\eta\neq 0}\frac{|\langle\xi,\eta\rangle|}{\|\eta\|_V}
=\sup_{\|\eta\|=1}\frac{|\langle\xi,\eta\rangle|}{\|\eta\|_V}, 
\end{align*}
the last equality holding because the function 
$$
\eta\in\{\eta\in H: \eta\neq 0\}\mapsto \frac{|\langle \xi,\eta\rangle|}{\|\eta\|_V}
$$
is homogeneous of degree zero.   
After taking reciprocals, we obtain 
\begin{align}\label{dnEq4}
r(V)&=
\inf_{\|\eta\|=1}\frac{\|\eta\|_V}{|\langle\xi,\eta\rangle|}.   
\end{align}
Now (\ref{dnEq4}) implies that there is 
a sequence of unit vectors  $\eta_n$ such that 
\begin{equation}\label{dnEq5}
\lim_{n\to\infty}\frac{\|\eta_n\|_V}{|\langle\xi,\eta_n\rangle|}=r(V).  
\end{equation}
Since 
$$
\frac{\|\eta_n\|_V}{|\langle\xi,\eta_n\rangle|}\geq \|\eta_n\|_V\geq r(V),    \qquad n=1,2,\dots,
$$
it follows that $\langle\xi,\eta_n\rangle\neq 0$ for large $n$; moreover,  since 
the left side converges to $r(V)$  we must have 
$$
\lim_{n\to\infty}\|\eta_n\|_V=r(V), \qquad {\rm{and \ }}\lim_{n\to\infty}|\langle\xi,\eta_n\rangle|=1.  
$$
Since 
$\xi$ and $\eta_n$ are unit vectors for which $|\langle \xi,\eta_n\rangle|$ converges 
to $1$, there is a sequence $\lambda_n\in\mathbb C$, $|\lambda_n|=1$, 
such that $\lambda_n\langle\xi,\eta\rangle=\langle\lambda_n\cdot\xi,\eta_n\rangle$ is nonnegative 
and converges to $1$.  It follows that 
$$
\lim_{n\to\infty}\|\lambda_n\cdot\xi-\eta_n\|^2=\lim_{n\to\infty}2-2\Re\langle\lambda_n\cdot\xi,\eta\rangle=0,   
$$
hence $\bar\lambda_n\cdot\eta_n$ converges in norm to $\xi$.  
By continuity of the norm $\|\cdot\|_V$, 
$$
\|\xi\|_V=\lim_{n\to\infty}\|\bar\lambda_n\cdot\eta_n\|_V=
\lim_{n\to\infty}\|\eta_n\|_V=r(V), 
$$
and (i) follows. 
\end{proof}

\begin{cor}\label{dnCor1}
If $r(V)>0$ then  $\|\cdot\|^V$ restricts to a bounded  norm-continuous 
function on the unit sphere of $H$ with the property that for every unit 
vector $\xi$,  $\|\xi\|^V=1$ iff $\xi\in V$ and $\|\xi\|^V=r(V)^{-1}$ 
iff $\xi$ is maximal.  
\end{cor}

\part{Normal states and normal functionals on $\mathcal B(H)$.}

Let $H$ be a Hilbert space and let $V$ be a norm-closed subset 
of the unit sphere of $H$ that satisfies axioms V1 and V2.  We now 
introduce a numerical function of normal states of $\mathcal B(H)$
that faithfully measures ``generalized entanglement", 
and we develop its basic properties in general.  When the inner radius 
of $V$ is positive, this function is shown to be the restriction of a {\em norm} 
on the predual $\mathcal B(H)_*$ to the space of normal states, or equivalently, the 
restriction of a norm on the Banach space $\mathcal L^1(H)$ of trace class 
operators to the space of density operators.

\section{Generalized entanglement of states}\label{S:es}

Fix a Hilbert space $H$.  
The Banach space $\mathcal B(H)_*$ of normal linear functionals on $\mathcal B(H)$ 
identifies naturally with the dual of the \cstar\ $\mathcal K$ of 
compact operators on $H$, and we may speak of the 
{\em weak$^*$-topology} on $\mathcal B(H)_*$.  Similarly, $\mathcal B(H)$ 
identifies with the dual of $\mathcal B(H)_*$,  and we may speak of the 
weak$^*$-topology on $\mathcal B(H)$.  Thus, a net of normal 
functionals $\rho_n$ converges weak$^*$ to zero iff 
$$
\lim_{n\to\infty}\rho_n(K)=0,\qquad \forall\ K\in\mathcal K, 
$$
and a net of operators $A_n\in\mathcal B(H)$ converges weak$^*$ to zero 
iff 
$$
\lim_{n\to\infty}\rho(A_n)=0,\qquad \forall \ \rho\in\mathcal B(H)_*.  
$$
There is a natural involution $\rho\mapsto\rho^*$ defined on $\mathcal B(H)_*$ by 
$$
\rho^*(A)=\overline{\rho(A^*)}, \qquad A\in\mathcal B(H), 
$$
and we may speak of {\em self adjoint} 
normal functionals $\rho$.  
Of course, $\mathcal B(H)_*$ identifies naturally with the Banach 
$*$-algebra 
of trace class operators, but that fact is 
not particularly useful for our purposes.  

Our aim is to introduce a measure of ``generalized entanglement" 
for normal states.   
It will be convenient to define 
it more generally as a function (\ref{esEq01}) 
defined on the larger Banach space $\mathcal B(H)_*$.  
For every $X\in\mathcal B(H)$, define 
\begin{equation}\label{esEq000}
\|X\|_V=\sup_{\xi,\eta\in V}|\langle X\xi,\eta\rangle|.   
\end{equation}
Axiom V2 implies that $\|\cdot\|_V$ is a norm, and obviously 
$\|X\|_V\leq \|X\|$ and $\|X^*\|=\|X\|$ for every $X$.  Consider 
the \cstar\ $\mathcal A$ obtained from 
the compact operators $\mathcal K\subseteq \mathcal B(H)$ 
by adjoining the identity operator 
$$
\mathcal A=\{K+\lambda\cdot\mathbf 1: K\in\mathcal K, \ \lambda\in \mathbb C\}.  
$$
Operators in $\mathcal A$ serve as ``test operators" for our purposes.  
The $V$-ball  in $\mathcal A$ 
\begin{equation}\label{esEq00}
\mathcal B=\{X\in \mathcal A:  \|X\|_V\leq 1\}
\end{equation}
is a norm-closed convex subset of $\mathcal A$ 
that is stable under the $*$-operation, 
stable under multiplication by complex scalars 
of absolute value $1$, 
 and it   
contains the unit ball of $\mathcal A$.  
Thus we can define 
an extended real-valued function $E:\mathcal B(H)_*\to[0,+\infty]$ by
\begin{equation}\label{esEq01}
E(\rho)=\sup_{X\in \mathcal B}\Re\rho(X)=\sup_{X\in \mathcal B}|\rho(X)|, 
\qquad \rho\in\mathcal B(H)_*.  
\end{equation}

\begin{rem}[Self adjoint elements of $\mathcal B(H)_*$]\label{esRem1}
Note that if $\rho=\rho^*$ is self adjoint functional in $\mathcal B(H)_*$, then 
$E(\rho)$ can be defined somewhat differently in terms of self adjoint operators:
\begin{align*}
E(\rho)&=\sup\{\Re\rho(X): X^*=X\in\mathcal A, \ \|X\|_V\leq 1\}
\\
&=\sup\{|\rho(X)|: X^*=X\in\mathcal A, 
\ \|X\|_V\leq 1\}.  
\end{align*}
Indeed, every $Z\in \mathcal B$ has a cartesian decomposition 
$Z=X+iY$ where $X$ and $Y$ are 
self adjoint with $X=(Z+Z^*)/2$, and we have  
$$
\Re\rho(Z)=\frac{1}{2}(\rho(Z)+\overline{\rho(Z)})=\frac{1}{2}\rho(Z+Z^*)=
\rho(X), 
$$
where $X=X^*\in\mathcal B$.  After noting $|\rho(X)|=\max(\rho(X), \rho(-X))$, 
we obtain
$$
E(\rho)\leq \sup\{|\rho(X)|: X^*=X\in \mathcal K, \ \|X\|_V\leq 1\}.  
$$
The opposite inequality is obvious.  
\end{rem}

In general, $E(\rho)$ can achieve the value $+\infty$
(see Remark \ref{crRem2}).  
We first determine when the set $\mathcal B$ 
is bounded.

\begin{prop}\label{ccProp2}
Let $r(V)$ be the inner radius of $V$ and let 
$\mathcal B_0$ be the set of all positive rank-one operators in $\mathcal B$.  Then 
\begin{equation}\label{esEq02}
\sup_{X\in B}\|X\|=\sup_{X\in B_0}\|X\|=\frac{1}{r(V)^2}.     
\end{equation}
Consequently, for every normal linear functional $\rho\in\mathcal B(H)_*$,   
\begin{equation}\label{esEq03}
\|\rho\|\leq E(\rho)\leq r(V)^{-2}\cdot\|\rho\|.  
\end{equation}
\end{prop}

\begin{proof}
To prove (\ref{esEq02}), 
it suffices to show that for every positive number $M$,  the following are equivalent:
\begin{enumerate}
\item[(i)]  $\|X\|\leq M\cdot\|X\|_V$ for every rank one projection $X\in \mathcal K$.  
\item[(ii)]  $\|X\|\leq M\cdot \|X\|_V$ for every $X\in \mathcal B(H)$.  
\item[(iii)]  $M\geq r(V)^{-2}$.   
\end{enumerate}
Since the implication (ii)$\implies$(i) is trivial, it is enough 
to prove (i)$\implies$(iii) and (iii)$\implies$(ii).  

(i)$\implies$(iii): 
Choose a unit vector $\zeta\in H$ and let $X$ be the rank
one projection $X\xi=\langle \xi,\zeta\rangle \zeta$, $\xi\in H$.  Then 
(i) implies 
\begin{align*}
1&=\|X\|\leq M\cdot\sup\{|\langle X\xi,\eta\rangle|: \xi,\eta\in V\}
\\
&=M\cdot\sup\{|\langle \xi,\zeta\rangle|\cdot|\langle \zeta,\eta\rangle|: \xi,\eta\in V\} 
\\
&=M\cdot \sup\{|\langle \zeta,\xi\rangle|^2: \xi\in V\}, 
\end{align*}
from which it follows that 
$$
\sqrt M\cdot\sup_{\xi\in V}\Re \langle \zeta,\xi\rangle=
\sqrt M\cdot\sup_{\xi\in V}|\langle \zeta,\xi\rangle|\geq 1.  
$$
Let $K$ be the closed convex hull of $V$.  After multiplying 
through by $\|\zeta\|$ for more general nonzero vectors $\zeta\in H$, 
the preceding inequality implies 
$$
\sqrt M\cdot\sup_{\xi\in K}\Re \langle \zeta,\xi\rangle =
\sqrt M\cdot\sup_{\xi\in V}\Re \langle\zeta,\xi\rangle \geq 
\|\zeta\|=\sup_{\|\eta\|\leq 1}\Re\langle \zeta,\eta\rangle.  
$$
Since every bounded real-linear functional $f: H\to\mathbb R$ must 
have the 
form $f(\xi)=\Re\langle \zeta,\xi\rangle$ for some vector $\zeta\in H$, 
a standard separation theorem implies that the unit ball 
of $H$ is contained in $\sqrt M\cdot K$, 
namely the closed convex hull of $\sqrt M\cdot V$.  
Hence $r(V)\geq M^{-1/2}$.  

(iii)$\implies$(ii):  Fix $X\in\mathcal B(H)$ and let $\xi_0,\eta_0\in H$ 
satisfy $\|\xi_0\|\leq 1$, $\|\eta_0\|\leq 1$.   
By definition of $r(V)$, hypothesis 
(iii) implies that both $\xi_0$ and $\eta_0$ belong to the 
closed convex hull of $\sqrt M\cdot V$, and hence 
\begin{align*}
|\langle X\xi_0,\eta_0\rangle|&\leq \sup\{|\langle X\xi,\eta\rangle|: 
\xi,\eta\in \sqrt M\cdot V\}
\\
&=M\cdot\sup_{\xi,\eta\in V}|\langle X\xi,\eta\rangle|=M\cdot\|X\|_V.  
\end{align*}
After taking the supremum over $\xi_0,\eta_0$, we obtain 
$\|X\|\leq M\cdot\|X\|_V$.  

The estimates (\ref{esEq03}) follow immediately from (\ref{esEq2}). 
\end{proof}

The basic properties of the function $E$ are summarized as follows.

\begin{thm}\label{esThm1}
The function $E:\mathcal B(H)_*\to[0,+\infty]$ satisfies:   
\begin{enumerate}
\item[(i)] For all $\rho_1,\rho_2\in \mathcal B(H)_*$,  $E(\rho_1+\rho_2)\leq E(\rho_1)+E(\rho_2)$.  
\item[(ii)]  For every nonzero $\lambda\in\mathbb C$ and every $\rho\in \mathcal B(H)_*$, 
$E(\lambda\cdot\rho)=|\lambda|\cdot E(\rho)$.  
\item[(iii)] $E$ is lower semicontinuous relative to the weak$^*$ topology of $\mathcal B(H)_*$.  
\item[(iv)] If $r(V)>0$, then $E$ is a {\em norm} equivalent 
to the norm of $\mathcal B(H)_*$.  
\end{enumerate}
Moreover, letting $\Sigma$ be the set of all normal states of $\mathcal B(H)$, we have 
\begin{equation}\label{esEq1}
\sup_{\rho\in\Sigma} E(\rho)=\sup_{X\in B}\|X\|=\frac{1}{r(V)^2},    
\end{equation}
the right side being interpreted as $+\infty$ when $r(V)=0$.  
\end{thm}

\begin{proof} (i), (ii) and (iii) are immediate consequences of 
the definition (\ref{esEq01}) of $E$ after noting that a supremum of continuous real-valued 
functions is lower semicontinuous, and (iv) follows from (\ref{esEq03}).  

To prove (\ref{esEq1}), let $\mathcal B_1=\{X=X^*\in\mathcal A: \|X\|_V\leq 1\}$ 
be the set of self adjoint operators in $\mathcal B$.  Remark \ref{esRem1} implies 
that 
$$
\sup_{\rho\in\Sigma} E(\rho)=\sup_{\rho\in\Sigma}\sup_{X\in \mathcal B_1}\rho(X)
=\sup_{X\in \mathcal B_1}\sup_{\rho\in\Sigma} \rho(X).  
$$
Noting that $\mathcal B_1=-\mathcal B_1$ and that the norm of a self adjoint operator agrees with 
its numerical radius, the right 
side can be replaced with 
$$
\sup_{X\in B_1}\sup_{\rho\in\Sigma}|\rho(X)|=\sup_{X\in B_1}\|X\|.  
$$
Formula (\ref{esEq1}) now follows from (\ref{esEq02}) of Proposition \ref{ccProp2}.  
\end{proof}

We may conclude that when the inner radius is positive, $E(\cdot)$ is uniformly continuous on 
the unit ball of $\mathcal B(H)_*$:  

\begin{cor}\label{esCor1}
Assume that $r(V)>0$.  Then for $\rho,\sigma\in \mathcal B(H)_*$ we have 
\begin{equation}\label{esEq2}
|E(\rho)-E(\sigma)|\leq r(V)^{-2}\cdot \|\rho-\sigma\|.  
\end{equation}
\end{cor}

\begin{proof}  Theorem \ref{esThm1} (iv) implies that $E(\cdot)$ is a norm on $\mathcal B(H)_*$, 
hence 
$$
|E(\rho)-E(\sigma)|\leq E(\rho-\sigma)\leq r(V)^{-2}\|\rho-\sigma\|, 
$$
the second inequality following from (\ref{esEq03}).  
\end{proof}

\section{$V$-correlated states and faithfulness of $E$}\label{S:f}

Given two unit vectors $\xi,\eta\in H$, we will write $\omega_{\xi,\eta}$ for 
the linear functional 
$$
\omega_{\xi,\eta}(A)=\langle A\xi,\eta\rangle, \qquad A\in\mathcal B(H).  
$$
One has $\|\omega_{\xi,\eta}\|=\|\xi\|\cdot\|\eta\|=1$, 
and $\omega_{\xi,\eta}^*=\omega_{\eta,\xi}$.  We begin by 
recalling two definitions from the introduction. 

\begin{defn} A normal state $\rho$ of $\mathcal B(H)$ is said to be 
{\em $V$-correlated} if for every $\epsilon>0$, there is an $n=1,2,\dots$, a set of 
vectors $\xi_1,\dots,\xi_n\in V$ and a set of positive reals $t_1,\dots,t_n$ 
satisfying $t_1+\cdots+t_n=1$ such that 
$$
\|\rho-(t_1\,\omega_{\xi_1,\xi_1}+\cdots+t_n\,\omega_{\xi_n,\xi_n})\|\leq \epsilon.
$$  
A normal state $\rho$ that is not $V$-correlated is said to be {\em entangled}.  
\end{defn}

By (\ref{esEq03}), $E(\rho)\geq 1$ for every normal state $\rho$.  
The purpose of this section is to 
prove the following result that characterizes 
entangled states by the inequality $E(\rho)>1$.    We assume that $H$ is a 
perhaps infinite dimensional Hilbert space, that 
$V\subseteq \{\xi\in H: \|\xi\|=1\}$ 
satisfies hypotheses V1 and V2, but we make no assumption about 
the inner radius of $V$.  

\begin{thm}\label{fThm1}
A normal state $\rho$ of $\mathcal B(H)$ is $V$-correlated iff $E(\rho)=1$.  
\end{thm}

The proof of Theorem \ref{fThm1} requires 
some preparation that is conveniently formulated in terms of the state space of the 
unital \cstar\ 
$$
\mathcal A=\mathcal K+\mathbb C\cdot\mathbf 1=
\{K+\lambda\cdot\mathbf 1: K\in\mathcal K,\ \lambda\in \mathbb C\},   
$$ 
which of course reduces to $\mathcal B(H)$ when $H$ is finite dimensional.  
After working out these preliminaries, 
we will return to the proof of Theorem 
\ref{fThm1} later in the section.  
The state space of $\mathcal A$ is compact convex in its relative weak$^*$-topology, 
not to be confused with the various weak$^*$-topologies described in the 
previous section.  
We write $\Sigma_V$ for 
the set of all states $\rho$ of $\mathcal A$ that satisfy
\begin{equation}\label{fEq1}
|\rho(X)|\leq \|X\|_V=\sup_{\xi,\eta\in V}|\langle X\xi,\eta\rangle|, 
\qquad  X\in\mathcal A.  
\end{equation}

\begin{thm}\label{fThm2}
Every state of $\Sigma_V$ is a weak$^*$-limit of states 
of $\mathcal A$ of the form 
$$
t_1\cdot\omega_{\xi_1,\xi_1}\restriction_\mathcal A+\cdots+t_n\cdot\omega_{\xi_n,\xi_n}\restriction_\mathcal A
$$
where $n=1,2,\dots$, $\xi_1,\dots,\xi_n\in V$ and the $t_k$ are positive reals with 
sum $1$.   
\end{thm}

\begin{proof}  Since (\ref{fEq1}) exhibits $\Sigma_V$ as an intersection of 
weak$^*$-closed subsets of the state space of $\mathcal A$, it follows that $\Sigma_V$ 
is weak$^*$-compact as well as convex.  The Krein-Milman theorem implies that 
$\Sigma_V$ is the weak$^*$-closed convex hull of its extreme points, hence it 
suffices to show that for every {\em extreme point} $\rho$ of $\Sigma_V$, there 
is a net of vectors $\xi_n\in V$ such that 
\begin{equation}\label{fEq2}
\rho(X)=\lim_{n\to\infty}\langle X\xi_n,\xi_n\rangle,\qquad X\in\mathcal A.  
\end{equation}

To that end, consider the somewhat larger set $\Omega_V$ of all 
bounded linear functionals $\omega$ on $\mathcal A$ that satisfy 
\begin{equation}\label{fEq3}
|\omega(X)|\leq \sup_{\xi,\eta\in V}|\langle X\xi,\eta\rangle|=\|X\|_V,\qquad X\in\mathcal A.  
\end{equation}
Since $\|X\|_V\leq \|X\|$, $\Omega_V$ is contained in the unit ball of the 
dual of $\mathcal A$, and it is clearly clearly 
convex and weak$^*$-closed, hence compact.  
We claim that 
\begin{equation}\label{fEq4}
\Omega_V=\overline{\co}^{\rm{weak}^*}\{\omega_{\xi,\eta}\restriction_\mathcal A: \xi,\eta\in V\}, 
\end{equation}
$\co$ denoting the convex hull.  
Indeed, the inclusion $\supseteq$ is immediate from the definition of $\Omega_V$.  
For the inclusion $\subseteq$, choose an operator $X\in \mathcal A$ and 
a real number $\alpha$ such that $\Re \omega_{\xi,\eta}(X)=\Re\langle X\xi,\eta\rangle\leq \alpha$ 
for all $\xi,\eta\in V$.  By axiom V1, this implies that for fixed $\xi,\eta\in V$ we have 
$$
|\langle X\xi,\eta\rangle|=\sup_{|\lambda|=1}\Re\lambda\langle X\xi,\eta\rangle=
\sup_{|\lambda|=1}\Re\langle X\lambda\cdot\xi,\eta\rangle\leq 
\sup_{\xi,\eta\in V}\Re\langle X\xi,\eta\rangle\leq \alpha
$$
and after taking the supremum over $\xi,\eta$ on the left side 
we obtain $\|X\|_V\leq \alpha$.  It follows that for every $\omega\in \Omega_V$,  
$$
|\omega(X)|\leq \|X\|_V\leq \alpha
$$
and (\ref{fEq4}) now follows from a standard separation theorem.  

Now let $\rho$ be an extreme point of $\Sigma_V$.  Then $\rho\in \Omega_V$, 
and we claim that in fact, $\rho$ is 
an extreme point of $\Omega_V$.   Indeed, if $\omega_1,\omega_2\in \Omega_V$ and 
$0<t<1$ are such that $\rho=t\cdot\omega_1+(1-t)\cdot\omega_2$, then 
$$
1=\rho(\mathbf 1)=t\cdot\omega_1(\mathbf 1)+(1-t)\cdot\omega_2(\mathbf 1).  
$$
Since $|\omega_k(\mathbf 1)|\leq \|\omega_k\|\leq 1$ and $1$ is an extreme point 
of the closed unit disk, it follows that $\omega_1(\mathbf 1)=\omega_2(\mathbf 1)=1$.  
Since $\|\omega_k\|\leq 1=\omega_k(\mathbf 1)$, this implies that both $\omega_1$ 
and $\omega_2$ are states of $\mathcal A$, hence $\omega_k\in\Sigma_V$.  By 
extremality of $\rho$, we conclude that $\omega_1=\omega_2=\rho$, as asserted.  

Finally, since $\rho$ is an extreme point of $\Omega_V$ and $\Omega_V$ is 
given by (\ref{fEq4}), Milman's converse of the Krein-Milman theorem implies 
that there is a net of pairs $\xi_n,\eta_n\in V$ such that 
$\omega_{\xi_n,\eta_n}$ converges to $\rho$ in the weak$^*$ topology.  

It remains to show that we can choose $\eta_n=\xi_n$ for all $n$, and for 
that consider $\omega_{\xi_n,\eta_n}(\mathbf 1)=\langle \xi_n,\eta_n\rangle$, 
which converges to $\rho(\mathbf 1)=1$ as $n\to\infty$.  This implies that 
$$
\|\xi_n-\eta_n\|^2=2(1-\Re\langle\xi_n,\eta_n\rangle)\to 0,
$$  
as $n\to\infty$, so that $\|\omega_{\xi_n,\xi_n}-\omega_{\xi_n,\eta_n}\|\leq \|\xi_n-\eta_n\|\to0$ 
as $n\to \infty$.  Hence $\omega_{\xi_n,\xi_n}$ converges 
weak$^*$ to $\rho$, and the desired conclusion (\ref{fEq2}) follows.  
\end{proof}

\begin{proof}[Proof of Theorem \ref{fThm1}]  It is clear from 
the definition (\ref{esEq01}) that 
$E(\rho)\geq 1$ in general.  We claim first that 
$E(\rho)=1$ for every $V$-correlated normal state $\rho$.  Indeed, 
since $E(\cdot)$ is a convex function that is lower semicontinuous with 
respect to the norm topology on states, 
the set $\mathcal C$ of all normal states $\rho$ for which $E(\rho)\leq 1$ 
is norm closed and convex.  It contains every state of the form $\omega_{\xi,\xi}$ 
for $\xi\in V$ since for every $X\in\mathcal A$ we have 
$$
\omega_{\xi,\xi}(X)\leq |\langle X\xi,\xi\rangle|\leq \sup_{\eta,\zeta\in V}|\langle X\eta,\zeta\rangle|
=\|X\|_V
$$
so that $E(\omega_{\xi,\xi})\leq 1$.  Hence $\mathcal C$ contains every 
$V$-correlated state.  

Conversely, let $\rho$ be a normal state for which $E(\rho)=1$, or equivalently,
$$
|\rho(X)|\leq \|X\|_V, \qquad X\in \mathcal A.
$$  
Theorem \ref{fThm2} implies that there is a net of normal states 
$\rho_n$ of $\mathcal B(H)$, each of which is a finite convex combination 
of states of the form $\omega_{\xi,\xi}$ with $\xi\in V$, such that 
$$
\rho(X)=\lim_{n\to\infty}\rho_n(X), \qquad X\in \mathcal A,   
$$
and in particular 
$$
\rho(K)=\lim_{n\to\infty}\rho_n(K),\qquad K\in\mathcal K.  
$$
It is well known that if a net of normal states converges to a normal state 
pointwise on compact operators, then in fact $\|\rho-\rho_n\|\to 0$ as 
$n\to\infty$ (for example, see Lemma 2.9.10 of \cite{arvMono}).  We conclude 
from the latter that $\rho$ is $V$-correlated.  
\end{proof}

\begin{rem} In the special case where $H$ is a tensor product of 
Hilbert spaces $H=H_1\otimes H_2$ and 
$V=\{\xi_1\otimes \xi_2: \xi_k\in H_k, \ \|\xi_1\|=\|\xi_2\|=1\}$, 
Holevo, Shirokov and Werner showed \cite{hskSep} that when $H_1$ and 
$H_2$ are infinite dimensional, there are normal 
states that can be norm approximated by convex combinations 
of vector states of the form $\omega_{\xi,\xi}$, $\xi\in V$, 
but which cannot be written as a discrete infinite convex combination 
$$
\rho=\sum_{k=1}^\infty t_k\cdot\omega_{\xi_k, \xi_k}
$$
with $\xi_k\in V$ and with nonnegative numbers $t_k$ having sum $1$.  On the other hand, 
they also show that every such $\rho$ can be expressed as  
an integral 
\begin{equation}\label{fEq5}
\rho(X)=\int_S\langle X\xi,\xi\rangle\,d\mu(\xi), \qquad X\in\mathcal B(H)
\end{equation}
where $\mu$ is a probability measure on the Polish 
space 
$$
S=\{\xi=\eta\otimes\zeta: \|\eta\|=\|\zeta\|=1\}.
$$  
It seems likely that an integral representation 
like (\ref{fEq5}) 
should persist for $V$-correlated states in the more general setting of Theorem \ref{fThm1},
where of course $S$ is replaced with $V$ -- though we have not pursued that issue.  
\end{rem}

\section{Maximally entangled states}\label{S:me}

The entanglement of a normal state $\rho$ satisfies 
$1\leq E(\rho)\leq r(V)^{-2}$, and the 
minimally entangled states were characterized as the $V$-correlated states in Theorem \ref{fThm1}. 
In this section we 
discuss states at the opposite extreme.  

\begin{defn}\label{meDef1}
A normal state $\rho$ satisfying  $E(\rho)=r(V)^{-2}$ 
is said to be {\em maximally entangled}.  
\end{defn}

We now calculate the entanglement of (normal) pure 
states in general,  and we characterize the maximally entangled 
pure states in cases where the inner radius of $V$ is positive.

\begin{thm}\label{meThm1}
Let $V$ be a norm-closed subset of the unit sphere of $H$ satisfying V1 and V2, 
let $\xi$ be a unit vector in $H$ and let $\omega$ the corresponding 
vector state $\omega(X)=\langle X\xi,\xi\rangle$, $X\in\mathcal B(H)$.  
Then 
\begin{equation}\label{meEq0}
E(\omega)=(\|\xi\|^V)^2.   
\end{equation}

Assuming further that $r(V)>0$, 
then $\omega$ is maximally entangled iff $\xi$ is a maximal vector.  
More generally, let $\rho$ be an arbitrary maximally entangled normal state, and 
decompose $\rho$ into a perhaps infinite convex combination 
of vector states 
\begin{equation}\label{meEq1}
\rho(X)=t_1\cdot\omega_1+t_2\cdot\omega_2+\cdots
\end{equation}
where the $t_k$ are positive numbers with sum $1$ 
and each $\omega_k$ has the form 
$\omega_k(X)=\langle X\xi_k,\xi_k\rangle$, 
$X\in\mathcal B(H)$, with $\|\xi_k\|=1$.  Then each $\omega_k$ 
is maximally entangled.  
\end{thm}

The proof of Theorem \ref{meThm1} makes use of the following basic inequality: 

\begin{lem}\label{meLem1}
For every $\xi,\eta\in H$ and every  $A\in\mathcal B(H)$, 
\begin{equation}\label{meEq3}
|\langle A\xi,\eta\rangle|\leq\|A\|_V \|\xi\|^V\|\eta\|^V.  
\end{equation}
\end{lem}

\begin{proof}[Proof of Lemma \ref{meLem1}]
After rescaling both $\xi$ and $\eta$, it is enough 
to show that 
\begin{equation}\label{meEq4}
\|\xi\|^V\leq 1,\quad \|\eta\|^V\leq 1\implies |\langle A\xi,\eta\rangle|\leq \|A\|_V.
\end{equation}
To that end, assume first that 
$\xi,\eta\in V$.  Then  
$$
|\langle A\xi,\eta\rangle|\leq \sup_{\xi,\eta\in V}|\langle A\xi,\eta\rangle|
=\|A\|_V.  
$$
Since $\langle A\xi,\eta\rangle$ is sesquilinear in $\xi,\eta$, 
the same inequality $|\langle A\xi,\eta\rangle|\leq \|A\|_V$ persists if $\xi$  and 
$\eta$ are finite 
convex combinations of elements of $V$, and by passing to 
the norm closure, $|\langle A\xi,\eta\rangle|\leq \|A\|_V$ 
remains true if $\xi$ and $\eta$ belong to the closed convex hull of $V$.  By 
Lemma \ref{dmLem1}, the closed convex hull of $V$ is $\{\zeta\in H: \|\zeta\|^V\leq 1\}$, 
and (\ref{meEq4}) follows.  
\end{proof}

\begin{proof} Let $\xi\in H$ be a unit vector with associated 
vector state $\omega$ and let $\mathcal A=\mathcal K+\mathbb C\cdot\mathbf 1$.  
Then for every $A\in\mathcal A$ satisfying 
$\|A\|_V\leq 1$, (\ref{meEq3}) implies 
$$
|\omega(A)|=|\langle A\xi,\xi\rangle|\leq (\|\xi\|^V)^2, 
$$
and $E(\omega)\leq (\|\xi\|^V)^2$ follows from the definition (\ref{esEq01}) after 
taking the supremum over $A$.  

To prove the opposite inequality $E(\omega)\geq (\|\xi\|^V)^2$, consider  
$$
\|\xi\|^V=\sup_{\|\zeta\|_V=1}|\langle\xi,\zeta\rangle|.     
$$
Let $\zeta_n$ be a sequence of vectors in $H$ satisfying 
satisfying $\|\zeta_n\|_V=1$ for all $n=1,2,\dots$ and 
$|\langle \xi,\zeta_n\rangle|\uparrow\|\xi\|^V$ as $n\to\infty$.  
Consider the sequence of rank one operators $A_1,A_2,\dots$ 
defined by $A_n(\eta)=\langle \eta,\zeta_n\rangle\zeta_n$, 
$\eta\in H$, and note that $\|A_n\|_V= 1$.  
Indeed, we have 
\begin{align*}
\|A_n\|_V&=\sup_{\eta_1, \eta_2\in V}|\langle A_n\eta_1,\eta_2\rangle|
=\sup_{\eta_1,\eta_2\in V}|\langle \eta_1,\zeta_n\rangle\langle \zeta_n,\eta_2\rangle|
\\
&=(\sup_{\eta\in V}|\langle\zeta_n,\eta\rangle|)^2=\|\zeta_n\|_V^2=1.  
\end{align*}
So by (\ref{esEq01}), $E(\rho)\geq |\rho(A_n)|$ for every $n=1,2,\dots$.  
But since 
$$
\rho(A_n)=\langle A_n\xi,\xi\rangle=|\langle \xi,\zeta_n\rangle|^2\uparrow (\|\xi\|^V)^2
$$
as $n\to\infty$, it follows that $E(\rho)\geq (\|\xi\|^V)^2$.

For the second paragraph, assume that $r(V)>0$.  Theorem 
\ref{dnThmA} implies that $\|\xi\|^V=r(V)^{-1}$ iff $\xi$ is a maximal 
vector; and from (\ref{meEq0}) we conclude that $\omega$ is 
a maximally entangled state iff $\xi$ is a maximal 
vector.  

let $\rho$ be a maximally 
entangled state of the form (\ref{meEq1}).  By symmetry and 
since all the $t_k$ are positive, 
it suffices to show that $\omega_1$ is maximally entangled.  
For that, consider the normal state  
$$
\sigma=\frac{t_2}{1-t_1}\omega_2+\frac{t_3}{1-t_1}\omega_3+\cdots.  
$$
We have $\rho=t_1\cdot\omega_1+(1-t_1)\cdot\sigma$, and since 
$E$ is a convex function, 
$$
\frac{1}{r(V)^2}=E(\rho)\leq t_1 E(\omega_1)+(1-t_1)E(\sigma). 
$$
Since $E(\omega_1)$ and $E(\sigma)$ are both  $\leq r(V)^{-2}$, 
it follows that $E(\omega_1)=E(\sigma)=r(V)^{-2}$, 
hence $\omega_1$ is a maximally entangled pure state.  
\end{proof}

\begin{rem}[Infinitely entangled states]\label{crRem2}
Consider the case $H=H_1\otimes H_2$ with $V$ the set of decomposable 
unit vectors $\eta_1\otimes\eta_2$, with $\eta_k$ a unit vector 
in $H_k$, $k=1,2$.  
When $\dim H_1=\dim H_2=\infty$, infinitely entangled normal states exist.  
Indeed, Proposition \ref{vvProp2} below implies that there are 
unit vectors $\xi$ satisfying 
$\|\xi\|^V=+\infty$ in this case, and by Theorem \ref{meThm1}, 
such a $\xi$ gives rise to a vector 
state $\omega$ for which $E(\omega)=+\infty$.  
\end{rem}

\part{$N$-fold tensor products}

In the remaining sections we consider Hilbert spaces 
presented as $N$-fold tensor products 
$$
H=H_1\otimes \cdots\otimes H_N
$$
in which at most one of the factors $H_k$ is infinite-dimensional.  We can 
arrange that the dimensions $n_k=\dim H_k$ 
increase $n_1\leq \cdots \leq n_N$, so that $n_{N-1}<\infty$.  
The set $V$ of distinguished vectors 
is the set of all decomposable unit vectors 
$$
V=\{\xi_1\otimes\cdots\otimes \xi_N: \xi_k\in H_k, \ \|\xi_1\|=\cdots\|\xi_N\|=1\}.  
$$

The general results above imply that we will have a rather complete 
understanding of separable states and entanglement once we 
determine the inner radius of $V$, have an 
explicit description of the maximal vectors,  
and identify the entanglement norm of states.  In the 
remaining sections we 
present our 
progress in carrying out those calculations.  We calculate the 
vector norms $\|\cdot\|_V$ and $\|\cdot\|^V$ and the entanglement measuring 
norm $E$ of normal states in general.  In order to determine the maximal vectors 
one must first calculate the inner radius $r(V)$.  While we are unable to 
obtain an explicit formula in general, we do obtain such a formula under 
the assumption that $H_N$ is ``large" in the sense that 
$n_N\geq n_1\cdots n_{N-1}$ and we characterize maximal vectors 
as those unit vectors that purify the tracial state of 
the subalgebra $\mathcal B(H_1\otimes \cdots\otimes H_{N-1})\otimes\mathbf 1_{H_N}$.  Of course, a natural setting in which all of 
the results of this section are valid is that in 
which exactly one of the factors of $H_1\otimes\cdots\otimes H_N$ is 
infinite dimensional.

\section{Calculation of the vector norms $\|\cdot\|_V$ and $\|\cdot\|^V$}\label{S:vv}

\begin{rem}[Projective tensor products]\label{vvRem1}
We begin by reviewing 
the definition and universal property 
of the projective tensor product $E_1\hat\otimes\cdots\hat\otimes E_N$ of 
complex Banach spaces $E_1, \dots, E_N$.  We require these results only when 
at most one of $E_1,\dots, E_N$ is infinite dimensional and we confine 
the discussion to such cases, with the $E_k$ arranged so that their 
dimensions $n_k=\dim E_k$ weakly increase 
with $k$ and satisfy $n_{N-1}<\infty$.  Every vector 
$z$ of the algebraic tensor product of vector spaces $E_1\odot \cdots\odot E_N$ can be expressed as 
a sum of elementary tensors 
\begin{equation}\label{crEq00}
z=\sum_{k=1}^n x^k_1\otimes\cdots \otimes x^k_N, 
\end{equation}
in many ways, with $1\leq n\leq n_1n_2\cdots n_{N-1}$, $x^k_j\in E_j$, $j=1,\dots,N$.  
The projective norm (or greatest cross norm)  $\|z\|_\gamma$ is defined as 
$$
\|z\|_\gamma=\inf \sum_{k=1}^n\|x^k_1\|\,\|x^k_2\|\cdots\|x^k_N\|
$$
the infimum extended over all representations of $z$ of the form (\ref{crEq00}).  
It is a fact that the norm $\|\cdot\|_\gamma$ makes 
the algebraic tensor product into a Banach space - the projective tensor product - 
denoted $E_1\hat\otimes\cdots\hat \otimes E_N$.   
The projective norm is a cross norm 
($\|x_1\otimes \cdots\otimes x_N\|_\gamma=\|x_1\|\cdots\|x_N\|$) that dominates 
every cross norm on $E_1\odot\cdots\odot E_N$.  

It is characterized by the following universal property: 
For every Banach space $F$ and every 
bounded multilinear mapping $B: E_1\times \cdots\times E_N\to F$, 
there is a unique bounded linear operator $L: E_1\hat\otimes \cdots\hat\otimes E_N\to F$ 
with the property  $L(x_1\otimes \cdots\otimes x_N)=B(x_1,\dots,x_N)$ for all $x_j\in E_j$, 
$1\leq j\leq N$, and the norm of the linearizing operator $L$ is given by 
$$
\|L\|=\sup\{\|B(x_1,\dots,x_N)\|: \|x_j\|\leq 1, j=1,\dots, N\}.
$$ 
In particular, the norm of a linear functional 
$F: E_1\hat\otimes\cdots\hat\otimes E_N\to \mathbb C$ is 
\begin{equation}\label{vvEq0}
\|F\|=\sup_{\|x_1\|=\cdots=\|x_N\|=1}|F(x_1\otimes x_2\otimes\cdots \otimes x_N)|.  
\end{equation}
Moreover, every bounded linear functional $F$ on $E_1\hat\otimes\cdots\hat\otimes E_N$ 
can be written as a finite linear combination of decomposable functionals 
\begin{equation}\label{vvEq0.1}
F=\sum_{k=1}^{n_1n_2\cdots n_{N-1}} F^k_1\otimes \cdots\otimes F^k_N, 
\end{equation}
where for each $j=1,\dots, N$, 
$F^k_j$ is a bounded linear functional on $E_j$.  
\end{rem}

We now calculate the vector norms $\|\cdot\|_V$ and $\|\cdot\|^V$ for cases 
in which $V$ is the set of decomposable unit vectors in $N$-fold tensor products 
$H_1\otimes \cdots\otimes H_N$ 
where the dimensions $n_k=\dim H_k$ weakly increase with $n_{N-1}<\infty$.  
The space $H_N$ is allowed 
to be infinite dimensional.  

\begin{thm}\label{vvThm1}
For every $\xi\in H_1\otimes \cdots\otimes H_N$, let $F_\xi$ be the 
element of the dual of the projective tensor 
product $H_1\hat\otimes \cdots \hat \otimes H_N$ defined by 
$$
F_\xi(\eta_1\otimes\cdots\eta_N)=\langle \eta_1\otimes \cdots \otimes \eta_N, \xi\rangle.  
$$
Then the norms $\|\cdot\|_V$ and $\|\cdot\|^V$ are given by 
\begin{equation}\label{vvEq1}
\|\xi\|_V = \|F_\xi\|, \quad \|\xi\|^V = \|\xi\|_\gamma, \qquad \xi\in H_1\otimes\cdots\otimes H_N.  
\end{equation}
\end{thm}

\begin{proof}
The first formula of (\ref{vvEq1}) is an immediate consequence of the definition 
of $\|\xi\|_V$ and the formula (\ref{vvEq0}), since 
\begin{align*}
\|\xi\|_V&=\sup_{\|\eta_1\|=\cdots=\|\eta_N\|=1}|\langle \eta_1\otimes\cdots\otimes\eta_N,\xi\rangle|
\\
&=\sup_{\|\eta_1\|=\cdots=\|\eta_N\|=1}|F_\xi(\eta_1\otimes\cdots\otimes \eta_N)|=\|F_\xi\|.  
\end{align*}

For the second formula,  write 
$$
\|\xi\|^V=\sup_{\|\eta\|_V\leq 1}|\langle \xi,\eta\rangle|=
\sup_{\|\eta\|_V\leq 1}|F_\eta(\xi)|.  
$$
The formula just proved asserts that $\|\eta\|_V=\|F_\eta\|$, 
so the preceding formula can be written 
$$
\|\xi\|^V=\sup_{\|F_\eta\|\leq 1}|F_\eta(\xi)|\leq \|\xi\|_\gamma.  
$$
For the opposite inequality, we use 
the Hahn-Banach theorem to find a linear functional $F$ of norm $1$ in the 
dual of $H_1\hat\otimes\cdots\hat\otimes H_N$ such that $\|\xi\|_\gamma=F(\xi)$.  
By the Riesz lemma there is a unique vector $\eta\in H_1\otimes\cdots\otimes H_N$ 
such that $F(\zeta)=\langle \zeta,\eta\rangle$ for all $\zeta$, and in particular 
$\|\xi\|_\gamma=F(\xi)=\langle \xi,\eta\rangle = F_\eta(\xi)$.  By the 
first part of the proof we have $\|\eta\|_V=\|F_\eta\|=1$.  
Hence 
$$
\|\xi\|_\gamma =\langle \xi,\eta\rangle \leq\sup_{\|\eta\|_V\leq 1}|\langle \xi,\eta\rangle|=\|\xi\|^V, 
$$
and $\|\xi\|_\gamma=\|\xi\|^V$ follows.  
\end{proof}

The following observation implies that $r(V)$ can be zero and 
infinitely entangled vectors can exist.  While the physics literature contains 
examples of infinitely entangled states 
(e.g., see \cite{keylEtAl}), it seems worthwhile to present concrete examples of 
that phenomenon in  this context.

\begin{prop}\label{vvProp2}  Consider the case $N=2$,  and let $H=H_1\otimes H_2$ 
where $H_1$ and $H_2$ are both infinite dimensional.  Then there are vectors $\xi\in H$ 
satisfying $\|\xi\|=1$ and $\|\xi\|^V=+\infty$.  
\end{prop}

\begin{proof}
Let $\theta_1, \theta_2,\dots$ be positive numbers with sum $1$, such 
as $\theta_k=2^{-k}$, let $n_1, n_2,\dots$ be positive integers 
such that $\theta_kn_k\to\infty$ as $k\to\infty$, and let $e_1,e_2,\dots$ 
and $f_1,f_2,\dots$ be orthonormal sets in $H_1$ and $H_2$ respectively.  
Partition the positive integers into disjoint subsets $S_1, S_2,\dots$ 
such that $|S_k|=n_k$ for $k=1,2,\dots$.  For every $k=1,2,\dots$, 
let $\xi_k$ be the vector 
$$
\xi_k=\sum_{j\in S_k}e_j\otimes f_j.  
$$
Obviously, $\|\xi_k\|^2=|S_k|=n_k$, and we claim that $\|\xi_k\|_V=1$. 
Indeed, 
$$
\|\xi_k\|_V=\sup_{\|\eta\|=\|\zeta\|=1}|\langle \xi_k,\eta\otimes\zeta\rangle|
=\sup_{\|\eta\|=\|\zeta\|=1}|\sum_{j\in S_k}\langle e_j,\eta\rangle\langle f_j,\zeta\rangle|
=1,  
$$ 
where the last equality is achieved with unit vectors $\eta, \zeta$ 
of the form 
$$
\eta=n_k^{-1/2}\sum_{k\in S_k}e_j, \qquad \zeta=n_k^{-1/2}\sum_{j\in S_k}f_j.  
$$
The vectors $\xi_1, \xi_2, \dots$ are mutually orthogonal, so that 
$$
\xi=\sum_k\frac{\sqrt\theta_k}{\|\xi_k\|}\,\xi_k=\sum_k\frac{\sqrt\theta_k}{\sqrt n_k}\,\xi_k
$$ 
defines a unit vector in $H$.  

We claim that $\|\xi\|^V=+\infty$.  To see that, fix $k=1,2,\dots$ and 
use $\|\xi_k\|_V=1$ to write 
$$
\|\xi\|^V=\sup_{\|\eta\|_V=1}|\langle \xi,\eta\rangle|\geq |\langle \xi,\xi_k\rangle|
=\frac{\sqrt\theta_k}{\sqrt n_k}\,\|\xi_k\|^2=\sqrt{\theta_kn_k}.  
$$
By the choice of $n_k$ the right side is unbounded, hence $\|\xi\|^V=+\infty$.  
\end{proof}

\section{Calculation of the entanglement norm $E$}\label{S:ms}

Continuing in the context of the previous section, we now 
calculate the entanglement norm $E(\rho)$ of normal states 
$\rho$ on $\mathcal B(H_1\otimes \cdots\otimes H_N)$.  We write 
$\mathcal L^1(H)$ for the Banach space of trace class operators on a 
Hilbert space $H$, with trace norm 
$$
\|A\|=\tr|A|, \qquad A\in\mathcal L^1(H), 
$$
$|A|$ denoting the positive square root of $A^*A$.  
Every normal linear functional $\rho$ on $\mathcal B(H)$ has 
a density operator $A\in \mathcal L^1(H)$,  defined by 
$$
\rho(B)=\tr(AB),\qquad B\in\mathcal B(H), 
$$
and the identification of $\rho$ with its density operator $A$ is a linear isometry.  

\begin{thm}\label{msThm1}
Let $\rho$ be a normal state of $\mathcal B(H_1\otimes \cdots\otimes H_N)$ with 
density operator $A$, $\rho(X)=\tr(AX)$.  The entanglement of $\rho$ is 
given by  
\begin{equation}\label{msEq1}
E(\rho)=\|A\|_\gamma,   
\end{equation}
where $\|\cdot\|_\gamma$ is the greatest cross norm on the projective tensor product 
of Banach spaces $\mathcal L^1(H_1)\hat\otimes\cdots\hat \otimes \mathcal L^1(H_N)$.  
\end{thm}

Before giving the proof, we first calculate the norm $\|B\|_V$,  defined on 
operators $B\in\mathcal B(H_1\otimes\cdots\otimes H_N)$  as in (\ref{esEq000}), 
in the current setting in which $V$ is 
the set of decomposable unit vectors of $H_1\otimes\cdots\otimes H_N$: 

\begin{lem}\label{msLem2}
For every operator $B\in \mathcal B(H_1\otimes \cdots\otimes H_N)$, one has 
\begin{equation}\label{msEq2}
\|B\|_V=\sup\{|\tr(B(T_1\otimes\cdots\otimes T_N))|: T_k\in \mathcal L^1(H_k), \ \tr|T_k|\leq 1\}.  
\end{equation}
\end{lem}

\begin{proof}
In this case, the  definition (\ref{esEq000}) of the norm $\|B\|_V$ becomes 
$$
\|B\|_V=\sup|\langle B(\xi_1\otimes\cdots\otimes \xi_N,\eta_1\otimes\cdots \otimes \eta_N\rangle|
$$
the supremum extended over all pairs $\xi_k, \eta_k\in H_k$, $k=1,\dots, N$ that 
satisfy $\|\xi_k\|= \|\eta_k\|= 1$.  It follows that this formula can be written 
equivalently as 
\begin{equation}\label{msEq3}
\|B\|_V=\sup\{|\tr(B(T_1\otimes\cdots\otimes T_N))|\}
\end{equation}
the supremum extended over all rank one operators $T_k\in \mathcal B(H_k)$ having norm $1$.  It is well known that 
for every Hilbert space $H$, the unit ball of the Banach space $\mathcal L^1(H)$ of 
trace class operators is the 
closure (in the trace norm) of the set of convex combinations of rank one operators 
of norm at most $1$.  It follows that the formula (\ref{msEq3}) is equivalent 
to (\ref{msEq2}).  
\end{proof}

\begin{proof}[Proof of Theorem \ref{msThm1}]  We claim first that the bounded linear 
functionals on the projective tensor product 
$\mathcal L^1(H_1)\hat\otimes \cdots\hat\otimes \mathcal L^1(H_N)$ are precisely those 
of the form 
\begin{equation}\label{msEq4}
F_B(A)=\tr (AB),\qquad A\in\mathcal L^1(H_1)\hat\otimes\cdots\hat \otimes \mathcal L^1(H_N), 
\end{equation}
where $B$ is a  operator in $\mathcal B(H_1\otimes \cdots\otimes H_N)$.  
Indeed, for every operator $B\in \mathcal B(H_1\otimes \cdots\otimes H_N)$, the universal 
property of the projective cross norm implies that there is a unique bounded linear 
functional $F_B$ on $\mathcal L^1(H_1)\hat\otimes\cdots\hat\otimes\mathcal L^1(H_N)$ that satisfies 
$$
F_B(T_1\otimes\cdots\otimes T_N)=\tr((T_1\otimes\cdots\otimes T_N)B), \quad 
T_k\in\mathcal L^1(H_k), \ 1\leq k\leq N.  
$$
For the opposite inclusion, by (\ref{vvEq0.1}), every bounded linear functional $F$ on 
$\mathcal L^1(H_1)\hat\otimes\cdots\hat\otimes L^1(H_N)$ is a finite sum of the form 
$$
F=\sum_{j=1}^n F^1_j\otimes\cdots\otimes F^N_j
$$
where $F^k_j$ belongs to the dual of $\mathcal L^1(H_k)$, $1\leq k\leq N$.  Letting 
$B^k_j\in \mathcal B(H_k)$ be the operator defined by $F^k_j(T)=\tr(TB^k_j)$, one sees 
that the operator 
$$
B=\sum_{j=1}^n B^1_j\otimes\cdots\otimes B^N_j\in \mathcal B(H_1\otimes\cdots\otimes H_N)
$$
satisfies (\ref{msEq4}), and the claim is proved.  

Note too that by the universal property of projective tensor products, Lemma 
\ref{msLem2} implies that the norm of the linear functional $F_B$ associated with 
an operator $B$ as in  (\ref{msEq4}) is given by 
\begin{equation}\label{msEq5}
\|F_B\|=\|B\|_V.  
\end{equation}

Fixing $\rho(X)=\tr(AX)$ as above, the Hahn-Banach theorem, together with the 
preceding remarks, implies that 
$$
\|A\|_\gamma=\sup_{\|F_B\|\leq 1}|F_B(A)|=\sup_{\|F_B\|\leq 1}|\tr(AB)|.  
$$
Using (\ref{msEq5}), the right side becomes 
$$
\sup_{\|F_B\|\leq 1}|\tr(AB)|=\sup_{\|B\|_V\leq 1}|\tr(AB)|=\sup_{\|B\|_V\leq 1}|\rho(B)| = E(\rho), 
$$
and (\ref{msEq1}) is proved.  
\end{proof}

\section{Calculation of the inner radius}\label{ci}

Let $V$ be the set of all decomposable unit vectors 
in a tensor product $H=H_1\otimes \cdots\otimes H_N$ 
with weakly increasing dimensions $n_k=\dim H_k$ and 
$n_{N-1}<\infty$.  In this section we establish a universal lower 
bound on $r(V)$, we show that this lower bound 
is achieved when $n_N$ is sufficiently large, and we 
exhibit maximal vectors for those cases.  

\begin{thm}\label{ciThm1}
In general, the inner radius satisfies  
\begin{equation}\label{ciEq1}
r(V_N)\geq \frac{1}{\sqrt{n_1n_2\cdots n_{N-1}}}.  
\end{equation}
\end{thm}

\begin{proof}  By formula (\ref{ccEq3}) of Theorem \ref{ccThm1}, it suffices to show that 
for every unit vector $\xi\in H$, 
\begin{equation}\label{ciEq2}
\|\xi\|^V\leq \sqrt{n_1n_2\cdots n_{N-1}}.  
\end{equation}
Fix orthonormal bases 
\begin{equation}\label{ciEq2.1}
\{e^1_1,\dots,e^1_{n_1}\}, \dots \{e^{N-1}_1,\dots,e^{N-1}_{n_{N-1}}\}
\end{equation}
for $H_1,\dots,H_{N-1}$ respectively.  
Every unit vector $\xi\in H_1\otimes\cdots\otimes H_N$ can be 
decomposed uniquely into a sum 
\begin{equation}\label{ciEq3}
\xi=\sum_{i_1=1}^{n_1}\cdots\sum_{i_{N-1}=1}^{N-1}e^1_{i_1}\otimes\cdots\otimes e^{N-1}_{i_{N-1}}\otimes
\xi_{i_1,\dots,i_{N-1}}, 
\end{equation}
where $\{\xi_{i_1,\dots,i_{N-1}}\}$ is a set of vectors in $H_N$ satisfying 
$$
\sum_{i_1,\dots,i_{N-1}=1}^{n_1,\dots,n_{N-1}}\|\xi_{i_1,\dots,i_{N-1}}\|^2=1.  
$$
Indeed, $\xi_{i_1,\dots,i_{N-1}}$ is the vector of $H_N$ defined by 
$$
\langle \xi_{i_1,\dots,i_{N-1}},\zeta\rangle=\langle \xi,e^1_{i_1}\otimes\cdots\otimes
e^{N-1}_{i_{N-1}}\otimes\zeta\rangle, \qquad \zeta\in H_N.  
$$
By Theorem \ref{vvThm1}, $\|\cdot\|^V$ is a cross norm on the 
algebraic tensor product $H_1\odot\cdots\odot H_N$, so from (\ref{ciEq3}) and 
the Schwarz inequality, we conclude that  
\begin{align*}
\|\xi\|^V&\leq 
\sum_{i_1,\dots,i_{N-1}}\|e_{i_1}^1\otimes\cdots\otimes e^{N-1}_{i_{N-1}}\otimes \xi_{i_1,\dots,i_{N-1}}\|^V
=\sum_{i_1,\dots,i_{N-1}}\|\xi_{i_1,\dots,i_{N-1}}\|
\\
&\leq (\sum_{i_1,\dots,i_{N-1}} 1 )^{1/2}( \sum_{i_1,\dots,i_{N-1}}
\|\xi_{i_1,\dots,i_{N-1}}\|^2)^{1/2}=(n_1\dots n_{N-1})^{1/2}, 
\end{align*}
and (\ref{ciEq2}) follows.  
\end{proof}

Assume now that $n_N\geq n_1n_2\cdots n_{N-1}$, choose a set of orthonormal bases $\{e^1_{i_1}\}, 
\dots, \{e^{N-1}_{i_{N-1}}\}$ for 
$H_1,\dots, H_{N-1}$ as in (\ref{ciEq2.1}), let 
$$
\{f_{i_1,\dots,i_{N-1}}: 1\leq i_1\leq n_1, \dots, 1\leq i_{N-1}\leq n_{N-1}\}
$$ 
be 
an orthonormal set in $H_N$, and consider the unit vector $\xi\in H_1\otimes\cdots\otimes H_N$ 
defined by 
\begin{equation}\label{ciEq4}
\xi=\frac{1}{\sqrt{n_1\cdots n_{N-1}}}\sum_{i_1=1}^{n_1}\cdots\sum_{i_{N-1}=1}^{n_{N-1}}
e^1_{i_1}\otimes\cdots\otimes e^{N-1}_{i_{N-1}}\otimes f_{i_1,\dots,i_{N-1}}.  
\end{equation}

\begin{thm}\label{ciThm2}
For all cases in which $n_N\geq n_1\cdots n_{N-1}$, we have 
\begin{equation}\label{ciEq5}
r(V)=\frac{1}{\sqrt{n_1\cdots n_{N-1}}}, 
\end{equation}
and vectors of the form (\ref{ciEq4}) are maximal vectors.  
\end{thm}

\begin{proof}
Let $\xi$ be a unit vector of the form (\ref{ciEq4}).  We will show that 
\begin{equation}\label{ciEq6}
\|\xi\|^V=\sqrt{n_1\cdots n_{N-1}}.
\end{equation}
Once (\ref{ciEq6}) is established, formula (\ref{ccEq3}) of 
Theorem \ref{ccThm1} implies  that 
$$
r(V)^{-1}=\sup_{\|\xi\|=1}\|\xi\|^V\geq \sqrt{n_1\cdots n_{N-1}}, 
$$
so that $r(V)\leq (n_1\cdots n_{N-1})^{-1/2}$, and (\ref{ciEq5}) will follow after an application 
of Theorem \ref{ciThm1}.  At that point, (\ref{ciEq6}) makes the assertion 
$\|\xi\|^V=r(V)^{-1}$, and 
Theorem \ref{dnThmA} will imply 
that $\xi$ is maximal.  

Thus it suffices to establish (\ref{ciEq6}). 
Note first that by (\ref{ccEq3}) and (\ref{ciEq1}), 
\begin{equation}\label{ciEq7}
\|\xi\|^V\leq \frac{1}{r(V)}\leq \sqrt{n_1\cdots n_{N-1}}.  
\end{equation}
For the opposite inequality, Theorem \ref{vvThm1}  
implies 
that $\|\xi\|^V$ is the projective cross norm $\|\xi\|_\gamma$, 
and it suffices to  exhibit a linear functional 
$F$ of norm $1$ on the projective tensor product 
$H_1\hat \otimes \cdots\hat \otimes H_N$ such that 
\begin{equation}\label{ciEq8}
F(\xi)=\sqrt{n_1\cdots n_{N-1}}.  
\end{equation}
For that, consider the vector 
$$
\eta=(n_1\cdots n_{N-1})^{1/2}\cdot\xi=
\sum_{i_1=1}^{n_1}\cdots\sum_{i_{N-1}}^{n_{N-1}}e^1_{i_1}\otimes \cdots\otimes e^{N-1}_{i_{N-1}}\otimes 
f_{i_1,\dots, i_{N-1}},  
$$ 
and define 
 $F$ on $H_1\hat \otimes  \cdots\hat\otimes H_N$ by 
$F(\zeta)=\langle \zeta,\eta\rangle$.  By the universal property of the projective 
tensor product, the norm of $F$ is 
$$
\|F\|=\sup\{|F(v_1\otimes\cdots\otimes v_N)|: v_k\in H_k, \ \|v_k\|\leq 1\}.  
$$
Choosing $v_k\in H_k$, we have 
\begin{align*}
F(v_1\otimes \cdots\otimes v_k)&=\langle v_1\otimes\cdots\otimes v_k,\eta\rangle 
\\
&=\sum_{i_1=1}^{n_1}\cdots\sum_{i_{N-1}=1}^{n_{N-1}}
\langle v_1,e^1_{i_1}\rangle\cdots\langle v_{N-1},e^{N-1}_{i_{N-1}}\rangle\langle v_N,f_{i_1,\dots,i_{N-1}}\rangle
\\
&
=\langle v_N, 
\sum_{i_1,\dots,i_{N-1}}\langle e^1_{i_1},v_1\rangle \cdots
\langle e^{N-1}_{i_{N-1}},v_{N-1}\rangle f_{i_1,\dots,i_{N-1}}\rangle.   
\end{align*}
Using orthonormality of $\{f_{i_1,\dots,i_{N-1}}\}$, we can write 
\begin{align*}
\sup_{\|v_N\|\leq 1}|F(v_1\otimes \cdots \otimes v_N)|&=
\|\sum_{i_1,\dots,i_{N-1}}\langle e^1_{i_1},v_1\rangle \cdots
\langle e^{N-1}_{i_{N-1}},v_{N-1}\rangle f_{i_1,\dots,i_{N-1}}\|^2
\\
&=
\sum_{i_1,\dots,i_{N-1}}|\langle e^1_{i_1},v_1\rangle |^2\cdots |\langle e^{N-1}_{i_{N-1}},v_{N-1}\rangle|^2
\\
&=\|v_1\|^2\cdots \|v_{N-1}\|^2,    
\end{align*}
so that  $\|F\|=\sup\{\|v_1\|^2\cdots \|v_{N-1}\|^2: \|v_k\|\leq 1\}=1$.  
Applying this linear functional to $\xi$, we find that 
$$
F(\xi)=\langle \xi,\eta\rangle=\sqrt{n_1\cdots n_{N-1}}\cdot\|\xi\|^2=\sqrt{n_1\cdots n_{N-1}}
$$
and the desired inequality $\|\xi\|_\gamma\geq \sqrt{n_1\cdots n_{N-1}}$ is proved.  
\end{proof}

\section{Significance of the formula $r(V)=(n_1n_2\cdots n_{N-1})^{-1/2}$}\label{S:sc}

Theorem \ref{ciThm1} asserts that for $N$-fold tensor products $H=H_1\otimes \cdots \otimes H_N$ in 
which the dimensions $n_k=\dim H_k$ increase with $k$ and satisfy $n_{N-1}<\infty$, 
the inner radius of the set $V$ of decomposable vectors satisfies 
\begin{equation}\label{scEq1}
r(V)\geq \frac{1}{\sqrt{n_1n_2\cdots n_{N-1}}}.   
\end{equation}
We have also seen that for fixed  
$n_1\leq \dots\leq n_{N-1}<\infty$, equality holds in (\ref{scEq1}) 
when $n_N$ is sufficiently large (see Theorem \ref{ciThm2}).  

In this 
section we show that equality in (\ref{scEq1}) can be 
characterized in a way that is perhaps unexpected, 
in that $r(V)=(n_1\cdots n_{N-1})^{-1/2}$ iff the tracial state 
of $\mathcal B(H_1\otimes\cdots\otimes H_{N-1})$ can be extended to 
a pure state of $\mathcal B(H_1\otimes \cdots\otimes H_{N-1})\otimes \mathcal B(H_N)$.  
We also characterize that situation in terms of the size of $n_N$.  

\begin{thm}\label{scThm1}
Let $V$ be the decomposable unit vectors in a tensor product of 
finite dimensional Hilbert spaces 
$H=H_1\otimes \cdots\otimes H_N$, with 
$n_k=\dim H_k$ weakly increasing with $k$, 
consider the subfactor  
$$
\mathcal A=\mathcal B(H_1\otimes\cdots\otimes H_{N-1})
$$
of $\mathcal B(H_1\otimes \cdots\otimes H_N)$,  
and let $\tau$ be the tracial state of $\mathcal A$. 
The following assertions are equivalent:
\begin{enumerate}
\item[(i)]  Minimality of the inner radius:  
\begin{equation}\label{scEq1.1}
r(V)=(n_1\cdots n_{N-1})^{-1/2}.
\end{equation}
\item[(ii)]  Existence of purifications: There is a unit vector $\xi\in H_1\otimes \cdots\otimes H_N$ such that 
\begin{equation}\label{scEq2}
\tau(A)=\langle (A\otimes \mathbf 1_{H_N})\xi,\xi\rangle, \qquad A\in \mathcal A.  
\end{equation}
\item[(iii)] Lower limit on $\dim H_N$: $n_N\geq n_1n_2\cdots n_{N-1}$.  
\end{enumerate}
\end{thm}

The proof of 
Theorem \ref{scThm1} requires the following elementary result.  

\begin{lem}\label{hLem2}  Let $H$ and $K$ be finite dimensional Hilbert spaces
and let $\omega$ be a faithful state of $\mathcal B(H)$.  If there 
is a vector $\xi\in H\otimes K$ such that 
$$
\omega(A)=\langle (A\otimes\mathbf 1_K)\xi,\xi\rangle, \qquad A\in \mathcal B(H), 
$$
then $\dim K\geq \dim H$.  
\end{lem}

\begin{proof}
Let $\eta$ be a unit vector in $H\otimes H$ such that 
$$
\omega(A)=\langle(A\otimes \mathbf 1_H)\eta,\eta\rangle=\langle 
(\mathbf 1_H\otimes A)\eta,\eta\rangle, \qquad A\in \mathcal B(H).  
$$
For example, setting $n=\dim H$, let $\Omega$ be the density operator of $\omega$, 
with eigenvalue list $\lambda_1\geq\cdots\geq \lambda_n>0$ and corresponding 
eigenvectors $e_1,\dots , e_n$.  One can take 
$$
\eta=\sqrt{\lambda_1}\cdot e_1\otimes e_1+\cdots+\sqrt{\lambda_n}\cdot e_n\otimes e_n.  
$$
Since $\omega$ is a faithful state,  $\eta$ is a cyclic and 
separating vector for $\mathcal B(H)\otimes \mathbf 1_H$.  

For every $A\in \mathcal B(H)$ we have 
$$
\|(A\otimes \mathbf 1_H)\eta\|^2=\omega(A^*A)=\|(A\otimes\mathbf 1_K)\xi\|^2,  
$$
hence there is an isometry $U: H\otimes H=(\mathcal B(H)\otimes\mathbf 1_H)\eta\to H\otimes K$ 
satisfying
$$
U(A\otimes\mathbf 1_H)\eta=(A\otimes\mathbf 1_K)\xi,\qquad A\in \mathcal B(H).
$$
It follows that $\dim H\cdot\dim K=\dim (H\otimes K)\geq \dim (H\otimes H)=(\dim H)^2$, 
and $\dim K\geq \dim H$ follows after canceling $\dim H$.  
\end{proof}

\begin{proof}[Proof of Theorem \ref{scThm1}]  
The implication (ii)$\implies$(iii) follows after applying  
Lemma \ref{hLem2} to the case 
$H=H_1\otimes\cdots\otimes H_{N-1}$ and $K=H_N$, 
and  (iii)$\implies$(i) is an immediate 
consequence of Theorem \ref{ciThm2}.

(i)$\implies$(ii):  Since $H_1\otimes\cdots\otimes H_N$ is finite dimensional, 
maximal vectors exist.    We claim that for every maximal vector $\xi$, one has   
\begin{equation}\label{scEq3}
\langle (E_1\otimes E_2\otimes \cdots\otimes E_{N-1}\otimes\mathbf 1_{H_N})\xi,\xi\rangle
=\frac{1}{n_1n_2\cdots n_{N-1}}.  
\end{equation}

For the proof, choose a unit vector $e_k\in E_kH_k$, $1\leq k\leq N-1$ and 
consider the operator $U: H_N\to H_1\otimes\cdots\otimes H_N$ defined by 
$$
U\zeta=e_1\otimes\cdots\otimes e_{N-1}\otimes \zeta, \qquad \zeta\in H_N.  
$$
$U$ is a partial isometry whose range projection is  
$E_1\otimes \cdots\otimes E_{N-1}\otimes\mathbf 1_{H_N}$, 
and since $\langle UU^*\xi,\xi\rangle=\|U^*\xi\|^2$, 
(\ref{scEq3}) is equivalent to the assertion
\begin{equation}\label{scEq4}
\|U^*\xi\|=\frac{1}{\sqrt{n_1\cdots n_{N-1}}}.  
\end{equation}
We claim first that $\|U^*\xi\|\leq (n_1\cdots n_{N-1})^{-1/2}$.  
Indeed, for every unit vector $\zeta\in H_N$ we have 
\begin{align*}
|\langle U^*\xi,\zeta\rangle|&=|\langle \xi,U\zeta\rangle|=
|\langle \xi,e_1\otimes\cdots\otimes e_{N-1}\otimes \zeta\rangle|
\\
&
\leq\sup_{\|\eta_1\|=\cdots=\|\eta_N\|=1}|\langle\xi,\eta_1\otimes\cdots\otimes\eta_N\rangle|
=\|\xi\|_V 
\\
&= r(V)=\frac{1}{\sqrt{n_1\cdots n_{N-1}}}.  
\end{align*}
where the equality $\|\xi\|_V=r(V)$ follows from the characterization of 
maximal vectors of Theorem \ref{dnThmA}.  
The 
asserted inequality now follows after taking the supremum over $\|\zeta\|=1$.  

To prove (\ref{scEq4}), choose orthonormal bases 
$$
\{e^1_1,\dots,e^1_{n_1}\}, \dots \{e^{N-1}_1,\dots,e^{N-1}_{n_{N-1}}\}
$$ 
for 
$H_1,\dots, H_{N-1}$ respectively, such that 
$e^1_1=e_1,\dots, e^{N-1}_1=e_{N-1}$.  For every sequence 
of integers $i_1,\dots,i_{N-1}$, 
$1\leq i_k\leq n_k$, consider the operator 
$$
U_{i_1,\dots,i_{N-1}}: \zeta\in H_N\mapsto e^1_{i_1}\otimes\cdots\otimes e^{N-1}_{i_{N-1}}\otimes \zeta
\in H_1\otimes\cdots\otimes H_N.  
$$
The preceding argument implies $\|U^*_{i_1,\dots,i_{N-1}}\xi\|^2\leq (n_1\cdots n_{N-1})^{-1}$
for each $i_1,\dots, i_{n_{N-1}}$, hence 
\begin{equation}\label{scEq5}
\sum_{i_1=1}^{n_1}\cdots\sum_{i_{N-1}=1}^{n_{N-1}}\|U^*_{i_1,\dots,i_{N-1}}\xi\|^2 \leq 
\sum_{i_1=1}^{n_1}\cdots\sum_{i_{N-1}=1}^{n_{N-1}}\frac{1}{n_1\cdots n_{N-1}}=1.  
\end{equation}
For each $k=1,\dots, N-1$ and every $i=1,\dots,n_k$, 
let $E^k_i$ be the projection onto the subspace of $H_k$ spanned 
by $e^k_i$.  For each $i_1,\dots,i_{N-1}$, we have 
$$
\|U^*_{i_1,\dots,i_{N-1}}\xi\|^2=\langle 
(E^1_{i_1}\otimes \cdots\otimes E^{n_{N-1}}_{i_{N-1}}\otimes\mathbf 1_{H_N})\xi,\xi\rangle. 
$$
Since the projections occurring in the right side are mutually orthogonal and 
sum to the identity operator of $H_1\otimes\cdots\otimes H_N$, 
the left side of (\ref{scEq5}) is 
$$
\sum_{i_1=1}^{n_1}\cdots\sum_{i_{N-1}=1}^{n_{N-1}}\langle 
(E^1_{i_1}\otimes \cdots\otimes E^{N-1}_{i_{N-1}}\otimes \mathbf 1_{H_n})\xi,\xi\rangle 
=\|\xi\|^2=1.  
$$
It follows that the inequality of (\ref{scEq5}) is actually equality;  
and since each summand satisfies 
$\|U^*_{i_1,\dots,i_{N-1}}\xi\|^2\leq (n_1\cdots n_{N-1})^{-1}$, 
we must have equality throughout the summands.  
Formula  (\ref{scEq3}) follows.  

Let $\mathcal S$ be the set of all operators $A\in \mathcal A$ for which 
(\ref{scEq2}) holds.  Obviously, $\mathcal S$ is a linear space, and by 
(\ref{scEq3}),  every  tensor product $E_1\otimes\cdots\otimes E_{N-1}$ 
of rank one projections $E_k\in \mathcal B(H_k)$ belongs to $\mathcal S$.   For 
fixed $k$, the rank one projections in $\mathcal B(H_k)$ span $\mathcal B(H_k)$, 
so by multilinearity, $\mathcal S$ contains all operators of the form 
$A_1\otimes \cdots\otimes A_{N-1}$ with $A_k\in \mathcal B(H_k)$.  Since 
operators of the form $A_1\otimes \cdots\otimes A_{N-1}$ span $\mathcal A$ 
itself,  Theorem \ref{scThm1} follows.  
\end{proof}

\begin{rem}[Finite dimensionality]\label{scRem1}
Notice that the hypothesis $n_N<\infty$ was used only in the proof 
of (i)$\implies$(ii), and there only to ensure the existence of maximal 
vectors.  If maximal vectors are known to exist in a setting in which 
$n_N=\infty$, then the proof of (i)$\implies$(ii) applies verbatim.  Of 
course, whenever (iii) holds, maximal vectors exist by Theorem \ref{ciThm2}.  
\end{rem}

\section{Homogeneity and the case $n_N\geq n_1\cdots n_{N-1}$}\label{S:h}

Continuing under the hypotheses $n_1\leq \cdots \leq n_{N-1}<\infty$, 
we show in this section that when $n_N\geq n_1\cdots n_{N-1}$, 
the set of maximal vectors in 
$H_1\otimes \cdots\otimes H_N$ is acted upon {\em transitively} by the 
unitary group of $H_N$, and we draw out several 
consequences.

\begin{thm}\label{hThm1}
Assume that $n_N\geq n_1\cdots n_{N-1}$ and 
let $\xi_1$ and $\xi_2$ be two maximal vectors in $H_1\otimes \cdots\otimes H_N$.  
Then there is a unitary operator $U$ in $\mathcal B(H_N)$ such that 
\begin{equation}\label{hEq1}
\xi_2=(\mathbf 1_{H_1}\otimes\cdots\otimes\mathbf 1_{H_{N-1}}\otimes U)\xi_1.  
\end{equation}
Maximal vectors are characterized 
as the unit vectors $\xi\in H_1\otimes \cdots\otimes H_N$ that purify the 
tracial state of $\mathcal B(H_1\otimes\cdots\otimes H_{N-1})$ 
as in (\ref{scEq2}).  

Finally, the maximal vectors 
for $H_1\otimes \cdots\otimes H_N$ are simply the vectors of the form 
\begin{equation}\label{hEq2}
\xi=\frac{1}{\sqrt{n_1n_2\cdots n_{N-1}}}
(
e_1\otimes f_1+\cdots +e_{n_1\cdots n_{N-1}}\otimes f_{n_1\cdots n_{N-1}}), 
\end{equation}
where $\{e_k: 1\leq k\leq n_1\cdots n_{N-1}\}$ is an orthonormal basis 
for $H_1\otimes\cdots\otimes H_{N-1}$ and 
$\{f_k: 1\leq k\leq n_1\cdots n_{N-1}\}$ is an arbitrary orthonormal set in $H_N$.  
\end{thm}

We require the following elementary consequence of 
familiar methods associated with the GNS construction.  
We sketch the proof for completeness.  

\begin{lem}\label{hLem1} 
Let $\xi_1$,  $\xi_2$ be vectors in $H_1\otimes\cdots\otimes H_N$ such that 
\begin{equation}\label{hEq3}
\langle (A\otimes\mathbf 1_{H_N})\xi_1,\xi_1\rangle =
\langle (A\otimes\mathbf 1_{H_N})\xi_2,\xi_2\rangle
\end{equation}
for all  
$A\in\mathcal B(H_1\otimes\cdots\otimes H_{N-1})$. 
Then there is a unitary operator 
$U\in \mathcal B(H_N)$ such that 
\begin{equation}\label{hEq4}
(\mathbf 1_{H_1\otimes\cdots\otimes H_{N-1}}\otimes U)\xi_1=\xi_2.
\end{equation}
\end{lem}

\begin{proof}
Consider the following subalgebra $\mathcal B$ of $\mathcal B(H_1\otimes\cdots\otimes H_N)$ 
$$
\mathcal B=\mathcal B(H_1\otimes \cdots\otimes H_{N-1})\otimes \mathbf 1_{H_N}.  
$$ 
$\mathcal B$ is a finite dimensional factor 
isomorphic to $\mathcal B(H_1\otimes \cdots\otimes H_{N-1})$ whose commutant is 
$\mathbf 1_{H_1\otimes \cdots\otimes H_{N-1}}\otimes\mathcal B(H_N)$.

For $k=1,2$, consider the finite dimensional subspace $H^k$ of the tensor 
product $H_1\otimes\cdots\otimes H_N$ 
defined by $L_k=\{B\xi_k: B\in\mathcal B\}$.  Since 
$$\|B\xi_k\|^2=\langle B^*B\xi_k,\xi_k\rangle =\langle B^*B\xi_2,\xi_2\rangle=\|B\xi_2\|^2, 
\qquad k=1,2,\quad B\in\mathcal B, 
$$ 
there is a unique partial isometry 
$V$ in the commutant of $\mathcal B$ having initial space $L_1$, 
final space $L_2$, such that 
$$
VB\xi_1=B\xi_2,\qquad B\in\mathcal B;  
$$
and in particular, this operator satisfies $V\xi_1=\xi_2$.  

Since both spaces $L_k$ are invariant under $\mathcal B$, they are the 
ranges of projections in the commutant of $\mathcal B$, and therefore 
must have 
the form $L_k=H_1\otimes\cdots\otimes H_{N-1}\otimes K_k$, $k=1,2$, 
where $K_k$ is a finite dimensional subspace of $H_N$.  Moreover, 
since $V$ belongs to the commutant of $\mathcal B$, 
it has the form $V=\mathbf 1_{H_1\otimes\cdots\otimes H_{N-1}}\otimes U_0$ 
where $U_0$ is a partial isometry in $\mathcal B(H_N)$ having initial 
and final spaces $K_1$ and $K_2$ respectively.  Finally, since a finite rank 
partial isometry $U_0\in \mathcal B(H_N)$ can always be extended to a unitary 
operator $U\in \mathcal B(H_N)$, we obtain a unitary 
operator $U\in\mathcal B(H_N)$ with the property asserted in (\ref{hEq4}).  
\end{proof}

\begin{proof}[Proof of Theorem \ref{hThm1}] 
Choose an orthonormal set 
$$
\{f_{i_1,\dots,i_{N-1}}:1\leq i_1\leq n_1,\dots,1\leq i_{N-1}\leq n_{N-1}\}
$$ 
in $H_N$ and let $\xi$ be the 
vector of the form (\ref{ciEq4}).  
Theorem \ref{ciThm2} implies that $\xi$ is a maximal vector.    

Let $\xi^\prime$ be another maximal vector.  The proof of the 
implication (i)$\implies$(ii) of  Theorem \ref{scThm1} implies that 
$$
\langle A\xi_1,\xi_1\rangle =\langle A\xi_2,\xi_2\rangle =\tau(A), \qquad A\in 
\mathcal B(H_1\otimes\cdots\otimes H_{N-1})\otimes\mathbf 1_{H_N}, 
$$
(see Remark \ref{scRem1}), where $\tau$ is the tracial state.  Lemma 
\ref{hLem1} implies that there is a unitary operator $U\in \mathcal B(H_N)$ 
such that $\xi^\prime=(\mathbf 1_{H_1\otimes\cdots\otimes H_{N-1}}\otimes U)\xi$.   
Notice that this implies 
that $\xi^\prime$ also has the form (\ref{hEq2}), in which 
$\{f_{i_1,\dots,i_{N-1}}\}$ is replaced with $\{Uf_{i_1,\dots,i_{N-1}}\}$. 
It also shows 
that every maximal vector purifies the tracial state $\tau$.  

Another application of Lemma \ref{hLem1} 
shows that every vector $\eta$ in the tensor product 
$H_1\otimes \cdots\otimes H_N$ that purifies 
the tracial state $\tau$ above must have the form 
$\eta=(\mathbf 1_{H_1\otimes\cdots\otimes H_{N-1}}\otimes U)\xi$ where 
$\xi$ is the vector above, therefore $\eta$ is also a maximal vector 
of the form (\ref{ciEq4}).  Finally, since every vector of the 
apparently more general form (\ref{hEq2}) must purify 
the tracial state $\tau$ as above, it follows from Lemma \ref{hLem1}  
that there is a unitary operator 
$U\in\mathcal B(H_N)$ such that $\eta=(\mathbf 1_{H_1\otimes\cdots \otimes H_{N-1}}\otimes U)\xi$.
It follows that  
$\eta$ can be rewritten so that it has the form (\ref{ciEq4}), and is therefore maximal.  
\end{proof}

\begin{rem}[Stability of maximal vectors when $\dim H_N$ is large]\label{hRem1}
It is of interest to reformulate the above results as follows.  Let $H$, $K$ 
be finite dimensional Hilbert spaces such that $\dim H\leq \dim K$,
consider the bipartite tensor product $G=H\otimes K$ with the associated 
set 
$$
V=\{\xi\otimes \eta\in G: \xi\in H,\ \eta\in K, \ \|\xi\|=\|\eta\|=1\}
$$ 
of unit decomposable vectors.  Suppose we are given a further decomposition 
of $H$ into a tensor product $H=H_1\otimes \cdots \otimes H_r$, with the resulting 
set 
$$
\tilde V=\{\xi_1\otimes\cdots\xi_r\otimes\eta\in G: \xi_k\in H_k,\  \eta\in K, \  \|\xi_k\|=\|\eta\|=1\}
$$
of decomposable unit vectors in $G$.  Then the preceding results 
show that {\em the sets $V$ and $\tilde V$ give rise to the same set of maximal vectors, and 
their  inner radii satisfy $r(V)=r(\tilde V)$. }

That fact seems remarkable, given that the entanglement 
measuring norms $\|\cdot\|^V$ and $\|\cdot\|^{\tilde V}$ are different.  Indeed, 
recall from Theorem \ref{vvThm1} that the entanglement measuring norms for $V$ and 
$\tilde V$ are, respectively, the projective cross norms on the bipartite 
tensor product $(H\otimes H)\hat \otimes K$ 
and the tripartite tensor product $H\hat\otimes H\hat\otimes K$, respectively.  
To see that the norms are different, it suffices to exhibit a linear functional 
$F$ on the vector space $H\odot H\odot K$ with the property that its norm in 
the dual of $(H\otimes H)\hat\otimes K$ is $1$ but its norm in the dual of 
$H\hat\otimes H\hat \otimes K$ is $<1$.  To that end, choose a unit vector 
$e\in H\otimes H$ that does {\em not} decompose into a tensor product $e_1\otimes e_2$,  
let $f$ be an arbitrary unit vector in $K$, and set 
$$
F(\xi_1\otimes \xi_2\otimes\eta)=\langle\xi_1\otimes \xi_2,e\rangle\langle \eta,f\rangle, 
\qquad \xi_k\in H, \quad \eta\in K.  
$$
If one views $F$ as a linear functional in the dual of $(H\otimes H)\hat\otimes K$, 
then its norm is $\|e\|\cdot\|f\|=1$.  On the other hand, 
\begin{align*}
\sup_{\|\xi_1\|=\|\xi_2\|=\|\eta\|=1}|F(\xi_1\otimes \xi_2 \otimes\eta) |&=
\sup_{\|\xi_1\|=\|\xi_2\|=\|\eta\|=1}|\langle\xi_1\otimes\xi_2,e\rangle|\cdot|\langle\eta,f\rangle|
\\
&=
\sup_{\|\xi_1\|=\|\xi_2\|=1}|\langle\xi_1\otimes\xi_2,e\rangle|<1, 
\end{align*}
since $e$ is not a decomposable vector.  This implies that the norm of $F$ 
as an element of the dual of $H\hat\otimes H\hat \otimes K$ is $<1$, and 
we conclude that $\|\cdot\|^V\neq \|\cdot\|^{\tilde V}$.  

In a similar way, one can see that while the entanglement measuring 
function $E$ of states is different for the two sets $V$ and $\tilde V$, 
the set of {\em maximally} entangled states is the same for both 
sets $V$ and $\tilde V$.  
\end{rem}

\section{Remarks on the case  $n_N<n_1n_2\cdots n_{N-1}$}\label{S:cn}

In this section we continue the discussion 
of $N$-fold tensor products $H=H_1\otimes \cdots \otimes H_N$ 
with increasing dimensions $n_k=\dim H_k$, with $n_{N-1}<\infty$, 
and with $V$ the set of decomposable unit vectors.  We have discussed 
the cases in which $n_N\geq n_1\cdots n_{N-1}$ at some length, having calculated the 
inner radius of $V$ and having identified the maximal vectors.  The following 
result and its corollary address the remaining cases.  
The fact is that we have little information about the inner radius 
and the structure of maximal vectors in such cases that goes beyond 
the content of Corollary \ref{cnCor1}.  Perhaps there 
is no simple formula for $r(V)$ in general.  

\begin{thm}\label{cnThm1}
If $n_N<n_1\cdots n_{N-1}$, then 
$$
r(V)>\frac{1}{\sqrt{n_1n_2\cdots n_{N-1}}}.  
$$
\end{thm}

\begin{proof}
By Theorem \ref{ciThm1}, $r(V)\geq (n_1\cdots n_{N-1})^{-1/2}$, and we have 
to show that equality cannot hold.  But if equality held, then 
the hypothesis on $n_N$ timplies that $H_1\otimes\cdots\otimes H_N$ is finite dimensional, 
so that maximal vectors exist.  Every maximal vector $\xi$ satisfies the criteria 
of Theorem \ref{scThm1} (i), but item (iii) of Theorem \ref{scThm1} contradicts 
the hypothesis on $n_N$.  
\end{proof}

\begin{cor}\label{cnCor1}
The inner radius is given by $r(V)=(n_1\cdots n_{N-1})^{-1/2}$ if $n_N\geq n_1\cdots n_{N-1}$;  
otherwise, $r(V)>(n_1\cdots n_{N-1})^{-1/2}$.  
\end{cor}

\begin{rem}[Best constants for the projective norm of 
$H_1\hat\otimes \cdots\hat\otimes H_N$]\label{cnRem2}
It is of interest to reformulate the information about the inner radius 
given by Theorems \ref{ciThm2} and 
\ref{cnThm1} in purely Banach space terms.  Given finite dimensional Hilbert spaces 
$H_1,\dots, H_N$, let $\|\cdot\|_\gamma$ be the projective cross norm on 
the tensor product 
$H_1\otimes \cdots \otimes H_N$ and let $\|\cdot\|$ be its Hilbert space 
norm.   Then one has the following information about the best 
constant $c$ for which $\|\xi\|_\gamma \leq c\cdot \|\xi\|$ for all 
$\xi\in H_1\otimes\cdots\otimes H_N$: 
$$
c=\sup_{\|\xi\|=1}\|\xi\|_\gamma = \sqrt{n_1\cdots n_{N-1}}, \qquad {\rm if \ }\ n_N\geq n_1\cdots n_{N-1}, 
$$
and 
$$
c=\sup_{\|\xi\|=1}\|\xi\|_\gamma< \sqrt{n_1\cdots n_{N-1}}, \qquad {\rm if \ }\ n_N < n_1\cdots n_{N-1}.  
$$
Note too that the preceding results provide no further information about 
the constant $c$ 
in cases where $n_N<n_1\cdots n_{N-1}$, and the problem of developing 
sharper information is one of obvious significance for quantum information theory 
as well as for the local theory of Banach spaces.  For example, in the case 
of bipartite tensor products, 
it is shown in \cite{GorLew}   that the space $\mathcal B(H_1,H_2)$ (endowed with 
the operator norm) fails to have local unconditional structure if the dimensions 
of $H_1$ and $H_2$ are
large, with further developments in \cite{GorNote}.  
Also see \cite{GorJunge}, an important paper on the local theory and the
many connections with ideal norms.  
\end{rem}

\begin{rem}[qubit triplets]\label{cnRem1}
The simplest case of tripartite tensor products to which our results 
do not apply is the case in which $V$ is the set of 
unit vectors $f\otimes g\otimes h$ in $\mathbb C^2\otimes \mathbb C^2\otimes \mathbb C^2$.  
We have not attempted to calculate $r(V)$ or determine the maximal 
vectors for this example; and if one seeks to extend the 
preceding calculations into the cases $n_N<n_1\cdots n_{N-1}$, this would seem the natural place to begin.   
Notice that Corollary \ref{cnCor1} implies $r(V)>2$.  

In a more qualitative direction, one might seek asymptotic 
information about the behavior of $r(V_N)$ for large $N$, 
where $V_N$ is the set of decomposable unit vectors of $(\mathbb C^2)^{\otimes N}$.  
\end{rem}

\section{Summary of results for $N$-fold tensor products}\label{S:sr}

We have not interpreted the main abstract results for multipartite tensor 
products.  For the reader's convenience, we 
conclude by summarizing the results of 
Proposition \ref{ccProp2}, and Theorems \ref{dnThmA}, \ref{fThm1}, \ref{fThm2}, \ref{meThm1}, 
\ref{vvThm1}, \ref{msThm1}, \ref{scThm1} in more concrete terms 
for these special cases.   
Let $H_1,\dots,H_N$ be Hilbert spaces whose dimensions $n_k=\dim H_k$ are 
weakly increasing, with $n_{N-1}<\infty$.  
For brevity, we confine ourselves to the case in which 
$n_N\geq n_1\cdots n_{N-1}$ where our results are sharp; however 
some of the following statements remain valid in the remaining 
cases as well.  What is missing in the remaining cases $n_N<n_1\cdots n_{N-1}$ is 
that we have only rough knowledge of the inner radius (see Theorem 
\ref{cnThm1}), and correspondingly little information about the structure 
of maximal vectors.  Obviously, the existence of those gaps in 
what we know about multipartite entanglement calls for further research.  

Let $V$ be the decomposable 
unit vectors 
$\xi_1\otimes\cdots\otimes\xi_N$ in the tensor product of Hilbert spaces $H=H_1\otimes\cdots\otimes H_N$, 
in which $\xi_k\in H_k$, and $\|\xi_k\|=1$.   

\begin{thm}\label{srThm1} Let $\|\cdot\|$ be the ambient norm of 
$H=H_1\otimes\cdots\otimes H_N$ and let $\|\cdot\|_\gamma$ be the 
norm of the projective tensor product of Hilbert spaces $H_1\hat\otimes\cdots\hat\otimes H_N$.  
The restriction of $\|\cdot\|_\gamma$  
to the unit sphere of $H$ 
$$
S=\{\xi\in H: \|\xi\|=1\}
$$  
has these properties.   Its range is the interval 
$\|S\|_\gamma=[1,\sqrt{n_1\cdots n_{N-1}}\,]$.  For every $\xi\in S$ 
one has $\|\xi\|_\gamma=1$ iff $\xi\in V$ is a decomposable vector, 
and 
$$
\|\xi\|_\gamma=\sqrt{n_1\cdots n_{N-1}}
\iff \xi \rm{\ is\ maximal} \iff \xi  \rm{\ has\ the\ form\ } (\ref{hEq2}).
$$  
The maximal 
vectors are also characterized as the unit vectors $\xi\in H$ that purify 
the tracial state of $\mathcal B(H_1\otimes\cdots\otimes H_{N-1})$ in the 
sense of (\ref{scEq2}).  

Let $\|\cdot\|^\gamma$ be the norm of the projective tensor product of Banach 
spaces $\mathcal L^1(H_1)\hat\otimes\cdots\hat\otimes \mathcal L^1(H_N)$, and 
let $\mathcal D$ be the space of all density operators - positive operators 
in $\mathcal B(H)$ having trace $1$.  The range of $\|\cdot\|^\gamma$ on $\mathcal D$ is 
$$
\|\mathcal D\|^\gamma=[1, n_1\cdots n_{N-1}].  
$$
Let $A\in\mathcal D$ and let $\rho(X)=\tr(AX)$ be the corresponding normal state 
of $\mathcal B(H)$.  Then $\rho$ is separable $\iff \|A\|^\gamma=1$, and for every rank 
one density operator $A\eta=\langle \eta,\xi\rangle\xi$, $\eta\in H$, 
$\|A\|^\gamma=n_1\cdots n_{N-1}\iff$ $\xi$ is a maximal vector.  
If a mixed state $\rho$ is maximally entangled in 
the sense that its density operator $A$ satisfies $\|A\|^\gamma=n_1\cdots n_{N-1}$, then 
$A$ is a convex combination of rank one projections associated with maximal vectors.   
\end{thm}

In particular, the unique entanglement measuring norms for vectors and states are identified 
in these cases as $\|\xi\|^V=\|\xi\|_\gamma$ and $E(\rho)=\|A\|^\gamma$, respectively, where 
$A$ is the density operator of the state $\rho$.

\bibliographystyle{alpha}

\begin{thebibliography}{PGWP{\etalchar{+}}08}

\bibitem[Arv03]{arvMono}
W.~Arveson.
\newblock {\em Noncommutative Dynamics and ${E}$-semigroups}.
\newblock Monographs in Mathematics. Springer-Verlag, New York, 2003.

\bibitem[Arv07]{arvEnt1}
W.~Arveson.
\newblock The probability of entanglement.
\newblock {\em preprint}, pages 1--31, 2007.
\newblock ar{X}iv:0712.4163.

\bibitem[Arv08]{arvEnt2}
W.~Arveson.
\newblock Quantum channels that preserve entanglement.
\newblock {\em preprint}, pages 1--14, 2008.
\newblock ar{X}iv:0801.2531.

\bibitem[BNT02]{BNT}
R.~A. Bertlmann, H.~Narnhofer, and W.~Thirring.
\newblock A geometric picture of entanglement and {B}ell inequalities.
\newblock {\em Phys. Rev. A}, 66:032319, 2002.

\bibitem[GJ99]{GorJunge}
Y.~Gordon and M.~Junge.
\newblock Volume ratios in {$L_p$} spaces.
\newblock {\em Studia Math}, 136:147--182, 1999.

\bibitem[GL74]{GorLew}
Y.~Gordon and D.~R. Lewis.
\newblock Absolutely summing operators and local unconditional structures.
\newblock {\em Acta Math}, 133:27--48, 1974.

\bibitem[Gor81]{GorNote}
Y.~Gordon.
\newblock A note on the {$GL$} constant of {$E\otimes_\epsilon F$}.
\newblock {\em Israel J. Math}, 39:141--144, 1981.

\bibitem[GRW08]{guhneEtAl}
O.~G{\"u}hne, M.~Reimpell, and R.~Werner.
\newblock Lower bounds on entanglement measures from incomplete information.
\newblock {\em preprint}, 2008.
\newblock ar{X}iv:0802.1734 [quant-ph].

\bibitem[HGBL05]{hylEtAl}
P.~Hyllus, O.~G{\"u}hne, D.~Bru{\ss}, and M.~Lewenstein.
\newblock Relations between entanglement witnesses and {B}ell inequalities.
\newblock {\em Phys. Rev. A}, 72:012321, 2005.

\bibitem[HHHH07]{HorSurvey}
R.~Horodecki, P.~Horodecki, M.~Horodecki, and K.~Horodecki.
\newblock Quantum entanglement.
\newblock {\em preprint}, 2007.
\newblock arXiv:quant-ph/0702225v2.

\bibitem[HSW05]{hskSep}
A.~S. Holevo, M.~E. Shirokov, and R.~Werner.
\newblock Separability and entanglement-breaking in infinite dimensions.
\newblock {\em preprint}, pages 1--12, 2005.
\newblock arXiv:quant-ph/0504204v1.

\bibitem[KSW02]{keylEtAl}
M.~Keyl, D.~Schlingemann, and R.~Werner.
\newblock Infinitely entangled states.
\newblock {\em preprint}, 2002.
\newblock ar{X}iv:quant-ph/0212014.

\bibitem[Per96]{peresSep}
A.~Peres.
\newblock Separability criterion for density matrices.
\newblock {\em Phys. Rev. Lett.}, 77(8):1413--1415, 1996.

\bibitem[PGWP{\etalchar{+}}08]{pgEtAl}
D.~Perez-Garcia, M.~Wolf, C.~Palazuelos, I.~Villanueva, and M.~Junge.
\newblock Unbounded violation of tripartite {B}ell inequalities.
\newblock {\em Comm. Math. Phys.}, 279:455--486, 2008.

\bibitem[Rud00]{rudSep}
O.~Rudolph.
\newblock A separability criterion for density operators.
\newblock {\em J. Phys. A: Math. Gen.}, 33:3951--3955, 2000.

\bibitem[Rud01]{rudEn}
O.~Rudolph.
\newblock A new class of entanglement measures.
\newblock {\em J. Math. Phys.}, 42:2507--2512, 2001.

\bibitem[WG03]{weiGo}
T.-C. Wei and P.~Golbart.
\newblock Geometric measure of entanglement and applications to bipartite and
  multipartite quantum states.
\newblock {\em Phys. Rev. A}, 68:042307, 2003.

\bibitem[WG07]{wangGu}
X.~Wang and S.~Gu.
\newblock Negativity, entanglement witnesses and quantum phase transition in
  spin-1 {H}eisenberg chains.
\newblock {\em J. Phys. A: Math. Theor.}, 40:10759--10767, 2007.

\end{thebibliography}
\newcommand{\etalchar}[1]{$^{#1}$}
\newcommand{\noopsort}[1]{} \newcommand{\printfirst}[2]{#1}
  \newcommand{\singleletter}[1]{#1} \newcommand{\switchargs}[2]{#2#1}

\end{document}